\documentclass [a4paper,10pt]{article}
\usepackage [utf8] {inputenc}
\usepackage {geometry}
\usepackage {amsmath}
\usepackage {amsthm}
\usepackage {amssymb}
\usepackage {amsfonts}
\usepackage{graphicx}
\graphicspath { {/home/tomos/Pictures/} } 
\title {An average theorem for tuples of $k$-free numbers in arithmetic progressions}
\date {Tomos Parry}
\begin {document}
\maketitle
\begin {center}
$\hspace {1mm}$
\\
\section {- \hspace {2mm}Introduction}
\end {center}
The Barban-Davenport-Halberstam (BDH) Theorem and its refinement due to Hooley and Montgomery (HM)
are important theorems in analytic number theory since they suggest what one believes to be 
the correct order of magnitude for the error term in the Prime Number Theorem for Arithmetic Progressions.  The 
question of additive patterns in prime numbers is also a central problem, but a theorem 
of BDH type is out of reach - indeed it is not even known, at the time of writing, that there are infinitely many primes $p$ such 
that $p+2$ is also prime, never mind the Hardy-Littlewood Conjecture.
\\
\\ Let $k\geq 2$.  If for a given $n$ there is no prime $p$ for 
which $p^k|n$ then $n$ is said to be $k$-free.  An asymptotic formula of similar 
shape to that in the BDH-HM Theorem is known for the $k$-free numbers, the current state of knowledge attained and summarised 
by Vaughan in \cite {kfree}.  Crucially, the corresponding question 
on additive patterns in the $k$-frees \emph {is} accessible.  For given non-negative 
integers $0\leq h_1<\cdot \cdot \cdot <h_r$ we call $n$ a \emph {$k$-free $r$-tuple 
associated to $\mathbf h:=(h_1,...,h_r)$} if the $n+h_i$ are all $k$-free, and 
write $\mathcal R=\mathcal R(\mathbf h)$ for the set of all such $n$.  The asymptotic count
\[ \sum _{n\leq x\atop {n\in \mathcal R}}1=\mathfrak Sx+\mathcal O\left (x^{2/(k+1)+\epsilon }\right ),\]
for some $\mathfrak S=\mathfrak S_\mathbf h>0$, is easily established (see \cite {mirsky}) 
and restricting to arithmetic progressions isn't too much harder.  Indeed in \cite {twins} twins in arithmetic progressions 
were investigated and it was shown easily that
\[ \sum _{n\leq x\atop {n\equiv a(q)}}\mu _k(n)\mu _k(n+1)=\eta (q,a)x+\mathcal O\left (x^{2/(k+1)+\epsilon }\right ),\]
for some $\eta (q,a)>0$, and that
\[ V(x,Q):=\sum _{q\leq Q}\sum _{a=1}^q
\left |\sum _{n\leq x\atop {n\equiv a(q)}}\mu _k(n)\mu _k(n+1)-\eta (q,a)x\right |^2\ll 
Q^2\left (\frac {x}{Q}\right )^{2/k+\epsilon }+x^{4/(k+1)+\epsilon };\]
and in \cite {meng} the method of Vaughan (that in \cite {kfree}) is followed to show that
\[ V(x,Q)\ll Q^2\left (\frac {x}{Q}\right )^{1/k+\epsilon }+x^{1+2/k}\log Q+x^{3/2+1/2k+\epsilon };\]
these results are important because the same results for primes are out of reach.
\\
\\ As far as we can see, however, there is no recorded asymptotic formula for this variance of $k$-free twins.  In 
this paper we achieve this, indeed for general tuples.
\\
\newtheorem* {walcott}{Theorem}
\begin {walcott}\label {walcott}
Fix natural numbers $k\geq 2$ and $r\geq 1$, denote by $\mathcal K$ the set of $k$-free numbers, fix 
non-negative integers $0\leq h_1<\cdot \cdot \cdot <h_r$, and let
\begin {equation}\label {tuple}
\mathcal R=\{ n\in \mathbb N|n+h_i\in \mathcal K,\hspace {2mm}i=1,...,r\} 
\end {equation}
be the set of $k$-free $r$-tuples associated to $\mathbf h:=(h_1,...,h_r)$.  Let for $q,a\in \mathbb N$ and $x,Q\geq 1$
\begin {eqnarray}\label {etaultra}
\eta (q,a)=\sum _{d_1,...,d_r=1\atop {(d_i^k,d_j^k)|h_i-h_j\hspace {1mm}(1\leq i,j\leq r)
\atop {(q,d_i^k)|a+h_i\hspace {1mm}(1\leq i\leq r)}}}^\infty 
\frac {\mu (d_1)\cdot \cdot \cdot \mu (d_r)}{[q,d_1^k,...,d_r^k]},
\hspace {10mm}E_x(q,a)=\sum _{n\leq x\atop {n\in \mathcal R\atop {n\equiv a(q)}}}1-\eta (q,a)x\hspace {15mm}
\end {eqnarray}
and
\begin {equation}\label {vary}
V(x,Q)=\sum _{q\leq Q}\sum _{a=1}^qE_x(q,a)^2.
\end {equation}
Take $\mathfrak c=\mathfrak c(r)$ to be any number in $[1/2,1)$ for which we know
\[ \int _{\pm \infty }\frac {|\zeta (\sigma +it)|^rdt}{(1+|t|)^{3/2}}\]
converges absolutely for all $\sigma \geq \mathfrak c$.  For 
each prime $p$ write $R_p$ for the number of different residues represented by 
the $h_1,...,h_r$ modulo $p^k$.  If always $R_p<p^k$ then for $1\leq Q\leq x$ and $\epsilon >0$
\[ V(x,Q)=Q^2\left (\frac {x}{Q}\right )^{1/k}P\Big (\log (x/Q)\Big )
+\mathcal O_{k,r,\mathbf h, \epsilon }\left (Q^2\left (\frac {x}{Q}\right )^{\mathfrak c(r)/k+\epsilon }+x^{1+2/(k+1)+\epsilon }\right )\]
where $P=P(r,k,\mathbf h)$ is a polynomial of degree at most $r-1$. 
\end {walcott}
As already mentioned the only theorems in this direction are upper bound result for twins of squarefree 
numbers.  In this case we can take $\mathfrak c=1/2$ since it is contained 
in classical results that for $\sigma \geq 1/2$
\[ \int _T^{2T}|\zeta (\sigma +it)|^2dt\ll _\sigma T\log T\]
and so our theorem then says
\[ V(x,Q)=x^{1/2}Q^{3/2}P\Big (\log (x/Q)\Big )
+\mathcal O_\epsilon \left (x^{1/4}Q^{7/4}+x^{5/3+\epsilon }\right )\]
for some linear function $P$.  Of course if the $h_i$ cover a complete residue system 
modulo some $p^k$ then there are no $k$-free $r$-tuples.
\\
\\ In \cite {kfree} the evaluation of the variance of $k$-free numbers 
is translated into a binary additive problem in $k$-free numbers which 
can be tackled with the circle method (following the general method of \cite {vaughangeneral}).  Aside 
from the last stage of the proof 
we use the method laid out there.  The main difficulty when comparing with \cite {kfree} is 
that the Gauss sum associated to $k$-free $r$-tuples is less accessible than that associated to $k$-frees; we use 
the methods of \cite {extremal} to get hold of this object, although a direct argument is also possible
.  
\\
\\ The paper is structured as follows: In Section 2 we collect 
the elementary facts about the distribution of $k$-free $r$-tuples in arithmetic progressions; in Section 3 we discuss the Gauss 
sum; in Section 4 we do most of the circle method work; in Section 5 we obtain the necessary results for the application 
of Perron's formula to the quantity remaining after the circle method work; and in Section 6 we carry out the main argument, using 
the results of the previous sections.
\\
\\ Throughout we consider $k\geq 2,$ $r\geq 1$ and $0\leq h_1<\cdot \cdot \cdot <h_r$ as fixed and write $\mathbf h=(h_1,...,h_r)$.  The 
implied constants in the $\mathcal O$ symbol will always be understood to 
be dependent on $k,r,\mathbf h$ and $\epsilon $, and $\epsilon $ may be taken to be arbitrarily small at each of its 
occurences.  Often (but always with explicit mention) we will write statements such as 
\[ ``f(X)\ll g(X)\]
where the $\ll $ containts terms up to $X^\epsilon $" -  here we mean 
\[ ``f(X)\ll X^\epsilon g(X)".\]
Whenever $s,\sigma $ and $t$ appear in the same context we will 
always mean a complex number $s$ with real and imaginary 
parts $\sigma $ and $t$.  We will write statements that involve $r$-tuples using vectors and mean that that statement 
is to hold for each vector component.  For 
example, $\nu \equiv \mathbf d\hspace {1mm}(\mathbf q)$ would mean $\nu \equiv d_i\hspace {1mm}(q_i)$ for each $i=1,...,R$, where 
the $R,d_i,q_i$ would be clear from context. A 
sum $\Sigma '^q$ will mean that the summation variables are restricted to numbers coprime to $q$.  The $R$-fold
divisor estimate $d_R(n)\ll n^\epsilon $ is well known, as 
is the (General) Chinese Remainder Theorem which says
\[ n\equiv \mathbf a\text { mod }(\mathbf q)\]
has exactly one solution modulo $[q_1,...,q_R]$ if $(q_i,q_j)|a_i-a_j$ and has no solutions otherwise.  We will use both these facts 
frequently but often forget to mention where they come from.  A coprimality condition may often disappear from one line to 
the next with the introduction of the Möbius function; here we are using
\[ \sum _{d|n}\mu (d)=\left \{ \begin {array}{ll}1&\text { if }n=1\\ 0&\text { if }n\not =1.\end {array}\right .\]
\begin {center}
\section {- \hspace {2mm}Counting $k$-free numbers in arithmetic progressions}
\end {center}
Counting $k$-free numbers amounts to counting solutions of congruences modulo $k$-th powers.  The precision we need for $r$-tuples 
is contained in \cite {mirsky} but we reproduce the proof since we need a slightly different result to the one stated there.
\\ \newtheorem {et}{Lemma}[section]
\begin {et}\label {et}
(i) For any $R,D\in \mathbb N$ with $D$ being $k$-free we have for $Z>0$
\[ \sum _{d_1,...,d_R\atop {[d_1^k,...,d_R^k,D]\leq Z}}1\ll _RZ^\epsilon \left (\frac {Z}{D}\right )^{1/k}.\]
(ii) For any $R\in \mathbb N$, any distinct $a_1,...,a_R\in \mathbb N_0$, any 
$\delta \in [0,1/3]$ and any $Y>0$
\[ \sum _{d_1\cdot \cdot \cdot d_R>Y\atop {(d_i^k,d_j^k)|a_i-a_j}}[d_1,...,d_R]^{k(\delta -1)}
\ll _{R,\mathbf a}Y^{1+k(\delta -1)+\epsilon }.\]
(iii) For $t\geq 1$, $R,d_1,...,d_R\in \mathbb N$ and distinct $a_1,...,a_R\in \mathbb N_0$ denote 
by $\mathcal N_{\mathbf d;\mathbf a}(t)$ the number of 
solutions $n\leq t$ to the system $n\equiv -\mathbf a\hspace {1mm}(\mathbf d^k)$.  Then for $Y>0$ we have
\begin {equation}\label {mirsky}
\sum _{d_1\cdot \cdot \cdot d_R>Y}\mathcal N_{\mathbf d;\mathbf a}(t)
\ll _{R,\mathbf a}t^\epsilon \Big (tY^{1-k+\epsilon }+t^{2/(k+1)}\Big ).
\end {equation}
\end {et}
\begin {proof}
(i) We have
\[ \sum _{[d_1^k,...,d_R^k,D]\leq Z}1=\sum _{[n^k,D]\leq Z}\sum _{[d_1,...,d_R]=n}1\ll _RZ^\epsilon \sum _{[n^k,D]\leq Z}1
=Z^\epsilon \sum _{l|D}\sum _{n^k\leq Zl/D\atop {(n^k,D)=l}}1.\]
Write $l_0=\prod _{p|l}p$.  Then the inner sum above is
\[ 
\sum _{n^k\leq Zl/Dl_0^k\atop {(n^kl_0^k,D)=l}}1\leq 
\left (\frac {Z}{D}\right )^{1/k}\frac {l^{1/k}}{l_0}\leq \left (\frac {Z}{D}\right )^{1/k}\]
since $D$ and therefore $l$ is $k$-free, and the claim follows for $D\leq Z$.  If $D>Z$ the LHS of the sum in question is zero.  
\\
\\ (ii) It is straightforward to establish with induction that for any $d_1,...,d_n\in \mathbb N$
\[ [d_1,...,d_n]\geq \frac {\prod _id_i}{\prod _{i\not =j}(d_i,d_j)}.\]
Therefore for any $d_1,...,d_R$ with $(d_i^k,d_j^k)|a_i-a_j$ 
\[ \frac {1}{[d_1^k,...,d_R^k]}\ll _{R,\mathbf a}\frac {1}{d_1^k\cdot \cdot \cdot d_R^k}\]
and so (since $k(\delta -1)
<-1$)
\begin {eqnarray*}
\sum _{d_1\cdot \cdot \cdot d_R>Y\atop {(d_i^k,d_j^k)|a_i-a_j}}[d_1,...,d_R]^{k(\delta -1)}
&\ll _{R,\mathbf a}&\sum _{d_1\cdot \cdot \cdot d_R>Y}(d_1\cdot \cdot \cdot d_R)^{k(\delta -1)}
\ll \sum _{n>Y}n^{k(\delta -1)+\epsilon }\ll Y^{1+k(\delta -1)+\epsilon }.
\end {eqnarray*}
\\
\\ (iii) We prove the claim by induction on $R$.  Suppose $t$ is larger than all the $a_1,...,a_R$ since 
otherwise the LHS of the sum in question 
is 
\[ \leq \sum _{d_1,...,d_R}\sum _{n\leq a_i\text { (some }i)\atop {d_i^k|n+a_i\text { (all }i)}}1\ll _{R,\mathbf a}1\]
anyway.  We have
\begin {eqnarray}\label {gawniweld}
\sum _{d>Y}\mathcal N_{d;a}(t)=\sum _{d>Y}\sum _{n\leq t\atop {n\equiv -a(d^k)}}1
\leq \sum _{Y<d\leq (t+a)^{1/k}}\left (\frac {t}{d^k}+1\right )\ll _atY^{1-k}+t^{1/k}\hspace {10mm}
\end {eqnarray}
which is (stronger than) the result for $R=1$ so suppose now the result 
holds for some $R$.  Let $Z>0$ be a parameter.  We have 
\begin {eqnarray*}
&&\sum _{d_1\cdot \cdot \cdot d_{R+1}>Y\atop {d_1\cdot \cdot \cdot d_R>Z}}\mathcal N_{d_1,...,d_{R+1};a_{1},...,a_{R+1}}(t)
\\ &&\hspace {20mm}\leq \hspace {4mm}\sum _{d_1\cdot \cdot \cdot d_R>Z}
\sum _{n\leq t\atop {n\equiv -a_i(d_i^k)\atop {i=1,...,R}}}
\sum _{d_{R+1}\atop {n\equiv -a_{R+1}(d_{R+1}^k)}}1
\\ &&\hspace {20mm}\ll \hspace {4mm}(t+a_{R+1})^{\epsilon }\sum _{d_1\cdot \cdot \cdot d_R>Z}
\mathcal N_{d_1,...,d_R;a_1,...,a_R}(t)
\\ &&\hspace {20mm}\ll _{R,\mathbf a}\hspace {4mm}t^{\epsilon }\Big (tZ^{1-k+\epsilon }+t^{2/(k+1)}\Big )
\end {eqnarray*}
by assumption, and since the argument would obviously be the same if we had the summation 
condition $d_1\cdot \cdot \cdot d_{R+1}/d_i>Z$ for some $1\leq i\leq R$ instead of $i=R+1$ we deduce 
\begin {eqnarray}\label {mirsky1}
\sum _{d_1,...,d_{R+1}>Y\atop {d_1\cdot \cdot \cdot d_{R+1}/d_i>Z\atop {\text {for some $i$}}}}
\mathcal N_{d_1,...,d_{R+1},a_1,...,a_{R+1}}(t)
&\ll _{R,\mathbf a}&t^\epsilon \Big (tZ^{1-k+\epsilon }+t^{2/(k+1)}\Big ).\hspace {15mm}
\end {eqnarray}
On the other hand if always $d_1\cdot \cdot \cdot d_{R+1}/d_i\leq Z$ then we 
must have $d_1\cdot \cdot \cdot d_{R+1}\leq Z^{1+1/R}$ so that 
\begin {eqnarray}\label {mirsky2}
&&\sum _{d_1\cdot \cdot \cdot d_{R+1}>Y\atop {d_1\cdot \cdot \cdot d_{R+1}/d_i\leq Z\atop {\text {for all $i$}}}}
\mathcal N_{d_1,...,d_{R+1};a_1,...,a_{R+1}}(t)\notag 
\\ &&\hspace {30mm}\leq \hspace {4mm}\sum _{Y<d_1\cdot \cdot \cdot d_{R+1}\leq Z^{1+1/R}}
\mathcal N_{d_1,...,d_{R+1};a_1,...,a_{R+1}}(t)\notag 
\\ &&\hspace {30mm}\leq \hspace {4mm}\sum _{Y<d_1\cdot \cdot \cdot d_{R+1}\leq Z^{1+1/R}\atop {(d_i^k,d_j^k)|a_i-a_j}}
\left (\frac {t}{[d_1^k,...,d_{R+1}^k]}+1\right )\notag 
\\ &&\hspace {30mm}\ll _{R,\mathbf a}\hspace {4mm}\Big (tY^{1-k+\epsilon }+Z^{1+1/R+\epsilon }\Big )
\end {eqnarray}
from part (ii).  Together \eqref {mirsky1} and \eqref {mirsky2} imply, assuming $Z\leq t$,
\begin {eqnarray*}
\sum _{d_1\cdot \cdot \cdot d_{R+1}>Y}
\mathcal N_{d_1,...,d_{R+1};a_1,...,a_{R+1}}(t)&\ll _{R,\mathbf a}&t^{\epsilon }
\bigg (tY^{1-k+\epsilon }+Z^{1+1/R}+tZ^{1-k}+t^{2/(k+1)}\bigg )
\\ &\ll &t^{\epsilon }\bigg (tY^{1-k+\epsilon }+t^{(R+1)/(Rk+1)}+t^{2/(k+1)}\bigg )
\end {eqnarray*}
having chosen $Z=t^{R/(Rk+1)}$.  The second term being less than the third, this is the result for $R+1$.
\end {proof}
\newtheorem {el}[et]{Lemma}
\begin {el}\label {el}
Let $\mathcal R$ be as in \eqref {tuple}, let $\eta (q,a)$ and $E_t(q,a)$ be as 
in \eqref {etaultra}, and let $\theta =1/k$ and $\Delta =2/(k+1)$.
\\
\\ (i)\hspace {5mm}For $t\geq 1$ and $q,a\in \mathbb N$ 
\[ \sum _{n\leq t\atop {n\in \mathcal R\atop {n\equiv a(q)}}}1=\eta (q,a)t
+\mathcal O\left (t^{\Delta +\epsilon }\right ).\]
(ii)\hspace {5mm}For $t,\gamma \geq 1$
\begin {eqnarray*}\hspace {10mm}&&\sum _{q\leq \gamma }\sum _{\nu =1}^q|E_t(q,\nu )|^2
\ll \gamma ^{2-2\theta }t^{2\theta }+t^{2\Delta }+\gamma t^\Delta,
\\ &&\sum _{q\leq \gamma }\frac {1}{q}\sum _{\nu =1}^q|E_t(q,\nu )|^2\ll \gamma ^{1-2\theta }t^{2\theta }+t^{2\Delta },
\\ \hspace {10mm}&&\sum _{q\leq \gamma }\frac {1}{q^{1-\theta }}\sum _{\nu =1}^q|E_t(q,\nu )|^2
\ll \gamma ^{1-\theta }t^{2\theta }+t^{2\Delta }+\gamma ^\theta t^\Delta 
\\ \text {and }&&\sum _{q\leq \gamma }\frac {1}{q^{2-\theta }}\sum _{\nu =1}^q|E_t(q,\nu )|^2\ll t^{2\Delta };
\end {eqnarray*}
here the $\ll $ symbol may include\footnote {see the notation explained 
in the introduction} terms of size $t^\epsilon ,\gamma ^\epsilon $.
\end {el}
\begin {proof}
Recall that $\mathbf h=(h_1,...,h_r)$ is fixed from the start.  A sum $\Sigma ^*$ over variables $d_1,...,d_r$ will mean that for 
all $i,j$ we have $(d_i^k,d_j^k)|h_i-h_j$.  For natural numbers $d_1,...,d_r$ write $d^*=[d_1,...,d_r]$.  For given 
$q,a,d_1,...,d_r\in \mathbb N$ and $t>0$ write 
$\mathcal N^{q,a}_{\mathbf d;\mathbf h}(t)$ for the number of solutions $n\leq t$ to 
the system of congruences $n\equiv -\mathbf h\text { mod }(\mathbf d^k)$ and $n\equiv a\text { mod }(q)$, and 
write $\mathcal N_{\mathbf d;\mathbf h}(t)=\mathcal N^{1,1}_{\mathbf d;\mathbf h}(t)$, as in 
Lemma \ref {et}.  Since it is well-known that for any $N\in \mathbb N$
\[ \sum _{d^k|N}\mu (d)=\left \{ \begin {array}{ll}1&\text { if $N$ is $k$-free}\\ 0&\text { if not}\end {array}\right .\]
we see from \eqref {tuple} that
\begin {equation}\label {wesentlich}
\sum _{n\leq t\atop {n\in \mathcal R\atop {n\equiv a(q)}}}1=
\sideset {}{^*}\sum _{d_1,...,d_r\atop {d_i\leq (t+h_i)^{1/k}}}\mu (\mathbf d)\mathcal N_{\mathbf d;\mathbf h}^{q,a}(t);
\end {equation}
here we obviously write $\mu (\mathbf d)=\mu (d_1)\cdot \cdot \cdot \mu (d_r)$.  
\\
\\ (i) From \eqref {wesentlich} we have for a parameter $Y\leq t^{1/k}$ to be chosen
\begin {eqnarray}\label {radio}
\sum _{n\leq t\atop {n\in \mathcal R\atop {n\equiv a(q)}}}1&=&
\sideset {}{^*}\sum _{d_1\cdot \cdot \cdot d_r\leq Y}\mu (\mathbf d)
\mathcal N_{\mathbf d;\mathbf h}^{q,a}(t)
+\mathcal O\left (\hspace {1.5mm}\sideset {}{^*}\sum _{d_1\cdot \cdot \cdot d_r>Y}
\mathcal N_{\mathbf d;\mathbf h}(t)\right )\notag 
\\ &=&\sideset {}{^*}\sum _{d_1\cdot \cdot \cdot d_r\leq Y}\mu (\mathbf d)
\mathcal N_{\mathbf d;\mathbf h}^{q,a}(t)
+\mathcal O\bigg (t^\epsilon \Big (tY^{1-k+\epsilon }+t^\Delta \Big )\bigg )\hspace {15mm}
\end {eqnarray}
from Lemma \ref {et} (iii).  The main term here is 
\begin {eqnarray*}
&&\sideset {}{^*}\sum _{d_1\cdot \cdot \cdot d_r\leq Y\atop {(q,\mathbf d^k)|a+\mathbf h}}\mu (\mathbf d)
\left (\frac {t}{[q,d^{*k}]}+\mathcal O(1)\right )
\\ &&\hspace {15mm}=\hspace {4mm}t\sideset {}{^*}\sum _{d_1,...,d_r\atop {(q,\mathbf d^k)|a+\mathbf h}}
\frac {\mu (\mathbf d)}{[q,d^{*k}]}+\mathcal O\left (t\sideset {}{^*}\sum _{d_1\cdot \cdot \cdot d_r>Y}
\frac {1}{d^{*k}}+\sum _{d_1\cdot \cdot \cdot d_r\leq Y}1\right )
\\ &&\hspace {15mm}=\hspace {4mm}t\eta (q,a)+\mathcal O\left (tY^{1-k+\epsilon }+Y^{1+\epsilon }\right )
\end {eqnarray*}
from \eqref {etaultra} and Lemma \ref {et} (ii), so \eqref {radio} becomes
\begin {eqnarray*}
&&\sum _{n\leq t\atop {n\in \mathcal R\atop {n\equiv a(q)}}}1
=t\eta (q,a)+\mathcal O\bigg (tY^{1-k+\epsilon }+Y^{1+\epsilon }+t^\epsilon \Big (tY^{1-k+\epsilon }+t^\Delta \Big )\bigg )
\end {eqnarray*}
which gives (i) on choosing $Y=t^{1/k}$.
\\
\\ (ii) From \eqref {wesentlich} we have for a parameter $Y\leq t^{1/k}$ to be chosen
\begin {eqnarray}\label {radio2}
\sum _{n\leq t\atop {n\in \mathcal R\atop {n\equiv \nu (q)}}}1&=&
\sideset {}{^*}\sum _{d_1\cdot \cdot \cdot d_r\leq Y}\mu (\mathbf d)
\mathcal N^{q,\nu }_{\mathbf d;\mathbf h}(t)
+\mathcal O\left (\sideset {}{^*}\sum _{d_1\cdot \cdot \cdot d_r>Y\atop {(q,\mathbf d^k)|\nu +\mathbf h}}
\sum _{n\leq t\atop {n\equiv -\mathbf h(\mathbf d^k)\atop {n\equiv \nu (q)}}}1\right ).\hspace {10mm}
\end {eqnarray}
The main term here is 
\begin {eqnarray*}
\sideset {}{^*}\sum _{d_1\cdot \cdot \cdot d_r\leq Y\atop {(q,\mathbf d^k)|\nu +\mathbf h}}\mu (\mathbf d)
\left (\frac {t}{[q,d^{*k}]}+\mathcal O(1)\right )&=&t\sideset {}{^*}
\sum _{d_1,...,d_r\atop {(q,\mathbf d^k)|\nu +\mathbf h}}\frac {\mu (\mathbf d)}{[q,d^{*k}]}
+\mathcal O\left (t\sideset {}{^*}\sum _{d_1\cdot \cdot \cdot d_r>Y\atop {(q,\mathbf d^k)|\nu +\mathbf h}}
\frac {1}{[q,d^{*k}]}+\sum _{d_1\cdot \cdot \cdot d_r\leq Y}1\right )
\\ &=&t\eta (q,\nu )+\mathcal O\left (t\sideset {}{^*}\sum _{d_1\cdot \cdot \cdot d_r>Y\atop {\nu \equiv -\mathbf h
\left ((q,\mathbf d^k)\right )}}\frac {1}{[q,d^{*k}]}+Y^{1+\epsilon }\right )
\end {eqnarray*}
from \eqref {etaultra}, therefore \eqref {radio2} implies 
\begin {eqnarray}\label {gaga}
\sum _{n\leq t\atop {n\in \mathcal R\atop {n\equiv \nu (q)}}}1-
t\eta (q,\nu )&\ll &\sideset {}{^*}\sum _{d_1\cdot \cdot \cdot d_r>Y\atop {\nu \equiv -\mathbf h\left ((q,\mathbf d^k)\right )}}\left (
\sum _{n\leq t\atop {n\equiv -\mathbf h(\mathbf d^k)\atop {n\equiv \nu (q)}}}1
+\frac {t}{[q,d^{*k}]}\right )+Y^{1+\epsilon }\notag 
\\ &=:&\mathcal T_Y(q,\nu )+Y^{1+\epsilon }.
\end {eqnarray}
In general for general positive functions $f,g$ and $S:=\sum _{d}\left \{ f(d)+g(d)\right \} $ it is easy to establish 
that $S^2\ll \sum _{d,d'}\left \{ f(d)f(d')+g(d)g(d')\right \} $.  Therefore
\begin {eqnarray*}
\sum _{\nu =1}^q|\mathcal T_Y(q,\nu )|^2&\leq &
\sideset {}{^*}\sum _{d_1\cdot \cdot \cdot d_r>Y\atop {d_1'\cdot \cdot \cdot d_r'>Y}}
\sum _{\nu =1\atop {\nu \equiv -\mathbf h(q,\mathbf d^k)\atop {\nu \equiv -\mathbf h(q,{\mathbf d'}^k)}}}^q\left (
\sum _{n,n'\leq t\atop {n\equiv -\mathbf h(\mathbf d^k)
\atop {n'\equiv -\mathbf h({\mathbf d'}^k)\atop {n\equiv n'\equiv \nu (q)}}}}1+\frac {t^2}{[q,{d^*}^k][q,{d'}^{*k}]}\right ).
\end {eqnarray*}
and the congruence conditions in the $\nu $ sum amount to one congruence modulo 
\[ [(q,d_1^k),...,(q,d_r^k),(q,{d_1'}^k),...,(q,{d_r'}^k)]=\frac {(q,d^{*k})(q,{d'}^{*k})}{(q,d^{*k},{d'}^{*k})}\]
so that the whole $\nu $ sum is
\[ \leq \sum _{\nu =1}^q\sum _{n,n'\leq t\atop {n\equiv -\mathbf h(\mathbf d^k)
\atop {n'\equiv -\mathbf h({\mathbf d'}^k)\atop {n\equiv n'\equiv \nu (q)}}}}1
+\frac {t^2q(q,d^{*k},{d'}^{*k})}{[q,{d^*}^k][q,{d'}^{*k}](q,d^{*k})(q,{d'}^{*k})}
=\sum _{n,n'\leq t\atop {n\equiv -\mathbf h(\mathbf d^k)
\atop {n'\equiv -\mathbf h({\mathbf d'}^k)\atop {n\equiv n'(q)}}}}1+\frac {t^2(q,d^{*k},{d'}^{*k})}{qd^{*k}{d'}^{*k}}\]
and therefore
\begin {eqnarray*}
\sum _{\nu =1}^q|\mathcal T_Y(q,\nu )|^2&\leq &
\sideset {}{^*}\sum _{d_1\cdot \cdot \cdot d_r>Y\atop {d_1'\cdot \cdot \cdot d_r'>Y}}
\left (\sum _{n,n'\leq t\atop {n\equiv -\mathbf h(\mathbf d^k)
\atop {n'\equiv -\mathbf h({\mathbf d'}^k)\atop {n\equiv n'(q)}}}}1+\frac {t^2(q,d^{*k},{d'}^{*k})}{qd^{*k}{d'}^{*k}}\right ).
\end {eqnarray*}
Since for any $N\in \mathbb N$
\[ \sum _{q\leq \gamma }\frac {(q,N)}{q}\ll N^\epsilon (\log \gamma +1)\ll N^\epsilon \gamma ^\epsilon \]
we see, on separating the terms with $n=n'$ since for these no divisor estimate is applicable, that
\begin {eqnarray}\label {huh}
\sum _{q\leq \gamma }\sum _{\nu =1}^q|\mathcal T_Y(q,\nu )|^2&\ll &\sideset {}{^*}
\sum _{d_1\cdot \cdot \cdot d_r>Y\atop {d_1'\cdot \cdot \cdot d_r'>Y}}\left (t^\epsilon 
\sum _{n,n'\leq t\atop {n\equiv -\mathbf h(\mathbf d^k)\atop {n'\equiv -\mathbf h({\mathbf d'}^k)}}}1
+\gamma \sum _{n\leq t\atop {n\equiv -\mathbf h(\mathbf d^k)\atop {n\equiv -\mathbf h({\mathbf d'}^k)}}}1
+\frac {t^2(d^{*k},{d'}^{*k})^\epsilon \gamma ^\epsilon }{d^{*k}{d'}^{*k}}\right )\notag 
\\ &\ll &t^\epsilon \left (\left (\hspace {1.5mm}\sum _{d_1\cdot \cdot \cdot d_r>Y}
\mathcal N_{\mathbf d,\mathbf h}(t)\right )^2+\gamma \mathcal N_{\mathbf d,\mathbf h}(t)\right )
+t^2\gamma ^\epsilon \left (\hspace {1.5mm}\sum _{d_1\cdot \cdot \cdot d_r>Y}d^{*k(\epsilon -1)}\right )^2\notag 
\\ &\ll &t^{\epsilon }\Bigg (t^2Y^{2-2k+2\epsilon }+t^{2\Delta }
+\gamma \Big (tY^{1-k+\epsilon }+t^{\Delta }\Big )\Bigg )+t^2\gamma ^\epsilon Y^{2-2k+4\epsilon }
\end {eqnarray}
from Lemma \ref {et} (iii) and (ii); in the second term in the second line we summed first over $\mathbf d'$ and 
used $\Sigma _{d'_i|n+h_i}1\ll t^\epsilon $, valid for 
large $t$; if $t$ is not large the first claim to be proven is clear since the obvious bound $\eta (q,\nu )\ll q^{\epsilon -1}$ means the LHS 
is then 
\[ \sum _{q\leq \gamma }\sum _{\nu =1}^q\left |\sum _{n\leq t\atop {n\equiv \nu (q)}}1+t^2q^{2\epsilon -2}\right |^2\ll \gamma .\]
Putting \eqref {huh} in \eqref {gaga} we get, assuming $Y\leq t$,
\begin {eqnarray*}
&&\sum _{q\leq \gamma }\sum _{\nu =1}^q\left |\sum _{n\leq t\atop {n\in \mathcal R\atop {n\equiv \nu (q)}}}1-
t\eta (q,\nu )\right |^2
\\ &&\hspace {10mm}\ll \hspace {4mm}t^{\epsilon }\Big (t^2Y^{2-2k+2\epsilon }+t^{2\Delta }+\gamma (tY^{1-k+\epsilon }+t^\Delta )\Big )
+t^2\gamma ^\epsilon Y^{2-2k+4\epsilon }+Y^{2+2\epsilon }\sum _{q\leq \gamma }\sum _{\nu =1}^q1\notag 
\\ &&\hspace {10mm}\ll \hspace {4mm}t^2Y^{2-2k}+t^{2\Delta }+\gamma (tY^{1-k}+t^\Delta )+Y^2\gamma ^2
\end {eqnarray*}
the $t^\epsilon ,\gamma ^\epsilon $ terms going into the $\ll $ symbol again.  Choosing $Y=(t/\gamma )^{1/k}$ gives 
the first claim and the others follow from partial summation.
\end {proof}
\begin {center}
\section {- \hspace {2mm}Gauss sums}
\end {center}
In this section we collect from \cite {extremal} 
the results needed to study the Gauss sum associated to $k$-free $r$-tuples. The letter $\mathcal S$ will always denote a general sequence whilst, as in the 
statement of our theorem, $\mathcal K$ denotes the $k$-free numbers and $\mathcal R$ the $k$-free $r$-tuples.  
\\
\\ If a sequence $\mathcal S$ satisfies for fixed $q$ and $a$
\[ \sum _{n\leq x\atop {n\in \mathcal S\atop {n\equiv a(q)}}}1=xf_\mathcal S(q,a)+E_\mathcal S(x;q,a)\]
for some $f(q,a)$ and some
\[ E_\mathcal S(x;q,a)=o(x),\hspace {10mm}x\rightarrow \infty ;\]
we say that \emph {$\mathcal S$ satisfies Criterion D}.  We define 
the \emph {density} of $\mathcal S$ as $\rho _\mathcal S=f(1,1)$ and if this is non-zero we define
\[ g_\mathcal S(q,a)=\frac {f_\mathcal S(q,a)}{\rho _\mathcal S}.\]
The \emph {Gauss sum} of $\mathcal S$ is defined as
\[ G_\mathcal S(q,a)=\sum _{\nu =1}^qe\left (\frac {a\nu }{q}\right )g_\mathcal S(q,a).\]
These definitions are all on page 92 of \cite {extremal}.  From Lemma 2.3 of \cite {extremal} (page 101) we can consider the Gauss sum as a function 
on $\mathbb Q/\mathbb Z$ and write 
\[ G_\mathcal S(a/q)=G_\mathcal S(q,a).\]
The Gauss sum 
is crucial to the exponential sum approximation in the circle method application later; indeed sorting the $n$ into arithmetic progressions 
modulo $q$ we see that for any $\mathcal S$ with non-zero density 
we have for $t>0$ 
\[ \sum _{n\leq t\atop {n\in \mathcal S}}e\left (\frac {an}{q}\right )=
\rho _\mathcal SG_\mathcal S(a/q)t+\sum _{\nu =1}^qe\left (\frac {a\nu }{q}\right )E_\mathcal S(t;q,\nu ).\]
Let $E_t(q,a)$ be as in \eqref {etaultra}.  Of course 
Lemma \ref {el} says that $E_t(q,a)=o(t)$ for fixed $q$ and $a$ so we must have $E_\mathcal R(t;q,a)=E_t(q,a)$.  Therefore 
the above says for any $t>0$
\begin {equation}\label {exponentialsum2}
\sum _{n\leq t\atop {n\in \mathcal R}}e\left (\frac {an}{q}\right )=\rho _\mathcal RG_\mathcal R(a/q)t+
\sum _{\nu =1}^qe\left (\frac {a\nu }{q}\right )E_t(q,\nu ),\hspace {10mm}\text {if }\rho _\mathcal R>0.
\end {equation}
For $h\in \mathbb N$ define the $h$-\emph {shift} of a sequence $\mathcal S$ as the sequence
\[ \{ n\in \mathbb N|n+h\in \mathcal S\} ,\]
which also obviously satisfies Criterion D, and write $G_\mathcal S^h$ for its Gauss 
sum.  From (2.18) and (2.19) of \cite {extremal} (page 108) we have for any $q,a\in \mathbb N$
\begin {equation}\label {shift}
G_\mathcal S^h(a/q)=e\left (-\frac {ah}{q}\right )G_\mathcal S(a/q).
\end {equation}
The rest of this section is concerned with evaluating the Gauss sum of the $k$-free $r$-tuples.  The underlying principle is that the Gauss sum associated 
to an intersection of sequences can be expressed in terms of the Gauss sums of the individual sequences via a convolution.  Since 
\begin {equation}\label {shifty}
\mathcal R=\bigcap _{i=1}^r\mathcal K_{h_i}
\end {equation}
and since the Gauss sum for the $k$-free numbers, and so from \eqref {shift} also the Gauss sum of their shifts, is accessible, we 
can therefore handle the Gauss sum of the $r$-tuples.
\\
\newtheorem {kfree}{Lemma}[section]
\begin {kfree}\label {kfree}
For prime $p$ and $a,l\in \mathbb N$ with $p\not |\hspace {1mm}a$
\[ G_\mathcal K(a/p^l)=\frac {-1}{p^k-1}\left \{ \begin {array}{ll}1&\text { if }l\leq k
\\ 0&\text { if }l>k.\end {array}\right .\]
\end {kfree}
\begin {proof}
The $k$-free numbers are $\{ s\in \mathbb N|\text { for all primes $p$, we have }p^k\not |\hspace {1mm}s\} $ so this 
is Lemma 5.3 of \cite {extremal} (page 128).
\end {proof}
As in definition (4.20) of \cite {extremal} (page 125) define, for a prime $p$, the \emph {$p$-local 
Gauss sum} $G_\mathcal S^p:\mathbb Q/\mathbb Z\rightarrow \mathbb C$ of $\mathcal S$ through
\[ G_\mathcal S^p(a/q)=\left \{ \begin {array}{ll}G_\mathcal S(a/q)&\text { if $a/q$ in lowest form has denominator 
a non-negative power of $p$}
\\ 0&\text { if not.}\end {array}\right .\]
As on page 118 of \cite {extremal} define, for two sequences $\mathcal S$ and $\mathcal S'$ satisfying 
Criterion $D$, the \emph {convolution} of $G_{\mathcal S}$ 
and $G_{\mathcal S'}$ as the function $G_\mathcal S\star G_{\mathcal S'}:\mathbb Q/\mathbb Z\rightarrow \mathbb C$ given through
\[ \left (G_\mathcal S\star G_{\mathcal S'}\right )(a/q)=\sum _{b/r\in \mathbb Q/\mathbb Z}
G_{\mathcal S'}(b/r)G_{\mathcal S}(a/q-b/r);\]
this is absolutely convergent by Lemma 1.1 of \cite {extremal} (page 92) so that in particular it is commutative, and for shifts 
of $k$-free numbers it must also be associative since in that case all the summations are finite summations, in view of \eqref {shift} 
and Lemma \ref {kfree}.
\newtheorem {gsl}[kfree]{Lemma}
\begin {gsl}\label {gsl}
Take $R\in \mathbb N$.  For any $H_1,...,H_R\in \mathbb N$ write $G_i$ for the Gauss sum 
of the $k$-frees shifted by $H_i$.  Then for any prime $p$ and any $a,l\in \mathbb N$ with $p\not |\hspace {1mm}a$ we 
have, if the $H_i$ are distinct modulo $p^k$,
\[ \left (G_1^p\star \cdot \cdot \cdot \star G_R^p\right )(a/p^l)=
\frac {-1/p^k}{(1-1/p^k)^R}\sum _{n=1}^{R}e\left (-\frac {aH_n}{p^l}\right )
\left \{ \begin {array}{ll}1&\text { if }l\leq k
\\ 0&\text { if }l>k.\end {array}\right .\]
and
\[ \left (G_1^p\star \cdot \cdot \cdot \star G_R^p\right )(0)=\frac {1-R/p^k}{(1-1/p^k)^R}.\]
\end {gsl}
\begin {proof}
We prove the first 
claim only, the proof of the second being essentially no different.  The 
result is clearly valid for $R=1$ in view of \eqref {shift} and Lemma \ref {kfree}, so 
suppose the result is true for some $R\in \mathbb N$ and take arbitrary distinct $H_1,...,H_{R+1}\in \mathbb N$.  Write $G$ for 
the Gauss sum of the $k$-frees.  For $b,L\in \mathbb N$ we have 
from \eqref {shift} 
\[ G_{R+1}(b/p^L)=
e\left (-\frac {bH_{R+1}}{p^L}\right )G(b/p^L)\]
so that, 
writing $q=p^l$,
\begin {eqnarray*}
\left (G_1^p\star \cdot \cdot \cdot \star G_{R+1}^p\right )(a/q)&=&\sum _{L=0}^\infty 
\sideset {}{'}\sum _{b=1}^{p^L}G^p_{R+1}(b/p^L)
\left (G_{1}^p\star \cdot \cdot \cdot \star G_R^p\right )\left (a/q-b/p^L\right )
\\ &=&\sum _{b=0}^{p^k}e\left (-\frac {bH_{R+1}}{p^k}\right )G^p(b/p^k)
\left (G_{1}^p\star \cdot \cdot \cdot \star G_R^p\right )\left (a/q-b/p^k\right ),
\end {eqnarray*}
the terms with $L>k$ vanishing in view of Lemma \ref {kfree}  
We now use the inductive 
hypothesis and (continue using) Lemma \ref {kfree}; we also drop the $p$ superscripts and write $g_n=(1-1/p^k)^n$.  Since 
clearly $G(0)=1$ the term $b=p^k$ contributes
\begin {eqnarray*}
G(0)\left (G_1\star \cdot \cdot \cdot \star G_R\right )(a/q)=\frac {-1/p^k}{g_R}\sum _{n=1}^{R}e\left (-\frac {aH_n}{q}\right ),
\end {eqnarray*}
the term $b=ap^k/q$ contribues
\begin {eqnarray*}
G(a/q)\left (G_1\star \cdot \cdot \cdot \star G_R\right )(0)
=\frac {-1/p^k}{g_1}\cdot \frac {1-R/p^k}{g_R}e\left (-\frac {aH_{R+1}}{q}\right )
\end {eqnarray*}
and the remaining terms contribute (both $G$ factors having non-trivial arguments)
\begin {eqnarray*}
&&\frac {-1/p^k}{g_1}\cdot 
\frac {-1/p^k}{g_R}\sum _{n=1}^{R}e\left (-\frac {aH_n}{q}\right )
\sum _{b=1\atop {b\not =ap^k/q,p^k}}^{p^k}e\left (-\frac {b(H_{R+1}-H_n)}{p^k}\right )
\\ &&\hspace {10mm}=\hspace {4mm}
\frac {1/p^{2k}}{g_{1+R}}\sum _{n=1}^{R}e\left (-\frac {aH_n}{q}\right )
\left (0-1-e\left (-\frac {a(H_{R+1}-H_n)}{q}\right )\right ),
\end {eqnarray*}
so the whole sum is
\begin {eqnarray*}
&&\left (\frac {-1/p^k}{g_R}-\frac {1/p^{2k}}{g_{1+R}}\right )\sum _{n=1}^{R}e\left (-\frac {aH_n}{q}\right )
+e\left (-\frac {aH_{R+1}}{q}\right )\left (\frac {-1/p^k}{g_1}\cdot \frac {1-R/p^k}{g_R}-\frac {1/p^{2k}}{g_{1+R}}
\sum _{n=1}^R1\right )
\\ &&\hspace {10mm}=\hspace {4mm}\frac {-g_1/p^k-1/p^{2k}}{g_{1+R}}\sum _{n=1}^Re\left (-\frac {aH_n}{q}\right )
+\frac {-1/p^k(1-R/p^k)-R/p^{2k}}{g_{1+R}}e\left (-\frac {aH_{R+1}}{q}\right )
\\ \\ &&\hspace {10mm}=\hspace {4mm}\frac {-1/p^k}{g_{1+R}}\sum _{n=1}^{R+1}e\left (-\frac {aH_n}{q}\right )
\end {eqnarray*}
and first claim follows.
\end {proof}
We now introduce the concept of \emph {quasi-multiplicativity} and introduce 
functions $G$ and $H$ which will be present throughout the paper, being essentially our exponential 
sum approximation for the $k$-free $r$-tuples.
\newtheorem {qm}[kfree]{Definition}
\begin {qm}\label {qm}
\emph {Take a function $f:\mathbb N^2\rightarrow \mathbb C$ for which $f(q,\cdot )$ has period $q$ for each $q\in \mathbb N$.  If for 
any $\omega \in \mathbb N$, any pairwise coprime $q_1,...,q_\omega \in \mathbb N$, and 
any $a_1,...,a_\omega \in \mathbb N$ we have
\[ f(q_1\cdot \cdot \cdot q_\omega ,a_1q/q_1+\cdot \cdot \cdot +a_\omega q/q_\omega )
=f(q_1,a_1)\cdot \cdot \cdot f(q_\omega ,a_\omega )\]
where $q=q_1\cdot \cdot \cdot q_\omega $, then we say that $f$ is \emph {quasi-multiplicative}; through induction 
this holds if and only if it holds for $\omega =2$.  We now take $q,a\in \mathbb N$ with $(a,q)=1$ and 
look at the value of $f(q,a)$ if $f$ is quasi-multiplicative.  Write $q=q_1\cdot \cdot \cdot q_\omega $ for 
the prime factorisation of $q$ and define $a_i$ through $a\equiv a_i$ mod $(q_i)$.  Then 
\[ a\equiv a_i\overline {q/q_1}q/q_1+\cdot \cdot \cdot +\overline {q/q_\omega }q/q_\omega \hspace {2mm}\text {mod }(q),\]
where $\overline {q/q_i}$ is inverse to $q/q_i$ mod $(q_i)$.  Therefore 
\[ f(q,a)=f\left (q_1,a_1\overline {q/q_1}\right )\cdot \cdot \cdot 
f\left (q_\omega ,a_\omega \overline {q/q_\omega }\right )\]
so that specifying the value of a quasi-multiplcative function $f(q,a)$ at prime powers $q$ and all $a$ with $(a,q)=1$ (and saying 
$f(1,1)=1$) is enough 
to determine $f$ for all $q\in \mathbb N$ and all $a$ with $(a,q)=1$.  Recall 
that the $h_1,...,h_r$ from our theorem are fixed from the outset.  For 
any prime $p$ denote by $H_1,...,H_R$ the $R=R_p$ different 
residues represented modulo $p^k$ by the $h_1,...,h_r$.  For any prime $p$ and any $a,l\in \mathbb N$ with $p\not |\hspace {1mm}a$ 
define
\[ G(p^l)=\frac {-1/p^k}{1-R/p^k}\left \{ \begin {array}{ll}
1&\text { if }l\leq k\\ 0&\text { if }l>k\end {array}\right .\]
and
\[ H(p^l,a)=\sum _{n=1}^{R}e\left (-\frac {aH_n}{p^l}\right )\left \{ \begin {array}{ll}
1&\text { if }l\leq k\\ 0&\text { if }l>k.\end {array}\right .\]
Define $G(q)$ and $H(q,a)$ for all $q,a\in \mathbb N$ with $(q,a)=1$ by extending 
multiplicatively and quasi-multiplcatively; note that $G$ is well-defined in view of the assumption in our 
theorem.$\hfill \blacksquare $}
\end {qm}
For large $p$ we have $R_p=r$ so $1-R_p/p^k\geq 1/2$ and so
\begin {equation}\label {babyleicht}
|G(p^t)|\leq \left \{ \begin {array}{ll}2/p^k&\text { for large $p$ }
\\ 0&\text { for all $p$ and }t>k,\end {array}\right .
\end {equation}
which we will use later, therefore for all $p$
\[ |G(p^t)|\ll \left \{ \begin {array}{ll}1/p^k&\text { always, in particular for $t\leq k$, }
\\ 0&\text { for }t>k.\end {array}\right .\]
We deduce $G(q)\ll 1/q$ for prime powers $q$ and so for general $q$
\begin {equation}\label {gbound}
G(q)\ll q^{\epsilon -1} ;
\end {equation}
also note $|H(q,a)|\leq R\ll 1$ holds for prime powers $q$ and for $a\in \mathbb N$ with $p\not |\hspace {1mm}$ so 
for general $q$ and $a$ with $(q,a)=1$
\begin {equation}\label {hbound}
H(q,a)\ll q^\epsilon .
\end {equation}
If a sequence $\mathcal S$ satisfying Criterion D has quasi-multiplicative $g_\mathcal S$ we say that \emph {$\mathcal S$ satisfies 
Criterion C}; see page 93 of \cite {extremal}; and we look 
at the intersection of such sequences.  In 
the paragraph containing equation (4.21) of \cite {extremal} (page 125) we have two sequences $\mathcal U$ and 
$\mathcal V$ satisfying Criterion $C$ with Gauss sums $\mathbf u$ and $\mathbf v$.  Shortly 
after $\mathbf w$ is defined as the Gauss sum of the intersection $\mathcal W:=\mathcal U\cap \mathcal V$ and then for any prime $p$
\[ \mathbf w_p=\frac {\mathbf u_p\star \mathbf v_p}{\left (\mathbf u_p\star \mathbf v_p\right )(0)}\]
according to (4.22) of \cite {extremal}, so long as $\rho _\mathcal W>0$
.  It follows for given $R\in \mathbb N$ that, if 
we have given sequences $\mathcal U_1,...,\mathcal U_R$ satisfying 
Criterion C with Gauss sums $\mathbf u_i$ and if $\mathbf w_R$ denotes the Gauss sum of the 
intersection $\mathcal W:=\mathcal S_1\cap \cdot \cdot \cdot \cap \mathcal S_R$, then for any prime $p$ and 
any $a/q\in \mathbb Q/\mathbb Z$ we have
\begin {eqnarray*}
\mathbf w_R^p(a/q)&=&\frac {\left (\mathbf u_1^p\star \cdot \cdot \cdot \star \mathbf u_R^p\right )(a/q)}
{\left (\mathbf u_1^p\star \cdot \cdot \cdot \star \mathbf u_R^p\right )(0)}
\end {eqnarray*}
so long as $\rho _\mathcal W>0$; moreover according to (1.4), (4.7) and Lemma 4.3 of \cite {extremal} (pages 92, 119 and 120) we have 
\[ \mathbf u_i\star \mathbf u_i=(\mathbf u_i\star \mathbf u_i)(0)\cdot \mathbf u_i
\]
and therefore we may drop repeated $p$-local Gauss sums from the above quotient of 
convolutions.  From \eqref {shift} and Lemma \ref {kfree} we 
have for any $h\in \mathbb N$, any prime $p$, any $l\geq 0$, and any $a\in \mathbb N$ with $p\not |\hspace {1mm}a$ 
\[ G_{\mathcal K_h}(a/p^l)=\left \{ \begin {array}{ll}e(-ah/p^l)G_{\mathcal K}(a/p^l)&\text { if }l\leq k
\\ 0&\text { if }l>k\end {array}\right \} =G_{\mathcal K_H}(a/p^l)\]
for any $H\in \mathbb N$ with $H\equiv h$ mod$(p^k)$. From \eqref {shifty} this 
discussion implies that for any prime $p$, any $l\geq 0$, and any $a\in \mathbb N$ with $p\not |\hspace {1mm}a$
\[ G_\mathcal R^p(a/p^l)=\frac {\left (G_{\mathcal K_{H_1}}^p\star \cdot \cdot \cdot \star G_{\mathcal K_{H_R}}^p\right )(a/p^l)}
{\left (G_{\mathcal K_{H_1}}^p\star \cdot \cdot \cdot \star G_{\mathcal K_{H_R}}^p\right )(0)}\]
where the $H_1,...,H_{R}$ are the $R=R_p$ different residues represented modulo $p^k$ by the $h_1,...,h_r$, and so 
from Lemma \ref {gsl} and Definition \ref {qm}
\begin {equation}\label {conv}
G_\mathcal R(p^l,a)
=\frac {-1/p^k}{1-R/p^k}\sum _{n=1}^{R}e\left (-\frac {aH_n}{q}\right )\left \{ \begin {array}{ll}
1&\text { if }l\leq k\\ 0&\text { if }l>k\end {array}\right \} =G(p^l)H(p^l,a)
\end {equation}
for any prime $p$, any $l\geq 0$, and any $a\in \mathbb N$ with $p\not |\hspace {1mm}a$; not to forget is that this is all 
subject to $\rho _\mathcal R>0$.  Moreover by Lemma 2.9 and Theorem 4.6 of \cite {extremal} (pages 110 and 125) it follows from 
\eqref {shifty} that $\mathcal R$ satisfies Criterion C and therefore, from 
Lemma 2.6 of \cite {extremal} (page 106), that $G_\mathcal R$ is quasi-multiplicative.  We deduce from Definition \ref {qm} that 
\[ G_\mathcal R(q,a)=G(q)H(q,a)\]
holds for $q,a\in \mathbb N$ with $(q,a)=1$ and so 
for general $q,a\in \mathbb N$
\begin {equation}\label {jajaja}
G_\mathcal R(q,a)=G\left (\frac {q}{(q,a)}\right )H\left (\frac {q}{(q,a)},\frac {a}{(q,a)}\right ),\hspace {10mm}\text {if }
\rho _\mathcal R>0. 
\end {equation}
We finish this section by establishing some easy properties of $G$ and $H$.
\newtheorem {cc}[kfree]{Lemma}
\begin {cc}\label {cc}
Define $G$ as in Definition \ref {qm} and write $\theta =1/k$.  For any $Z\geq 1$
\[ \sum _{Z<q\leq 2Z}|G(q)|\ll Z^{\theta -1+\epsilon }.\]
This implies in particular
\[ \sum _{q\leq Z}q|G(q)|\ll Z^{\theta +\epsilon },\hspace {5mm}\sum _{q\leq Z}q^2|G(q)|^2\ll Z^{\theta +\epsilon },\]
\[ \sum _{q>Z}|G(q)|\ll Z^{\theta -1+\epsilon },
\hspace {5mm}\sum _{q>Z}|G(q)|^2\ll Z^{\theta -2+\epsilon }\hspace {5mm}
\sum _{q>Z}q^{1+\epsilon }|G(q)|^2\ll Z^{\theta -1+\epsilon }\]
and 
\[ \sum _{q\leq Z}q^{2-\theta }|G(q)|^2\ll Z^{\epsilon }.\]
\end {cc}
\begin {proof}
By \eqref {babyleicht} we have 
\[ p^{t(1-1/k)}|G(p^t)|\leq \left \{ \begin {array}{ll}2/p&\text { for large $p$ and }1\leq t\leq k
\\ 0&\text { for all $p$ and }t>k\end {array}\right .\]
therefore by multiplicativity
\begin {eqnarray*}
\sum _{q\leq 2Z}q^{1-1/k}|G(q)|&\leq &
\prod _{p\leq 2Z}\left (1+\sum _{t\geq 1}p^{t(1-1/k)}|G(p^t)|\right )
\\ &\leq &\prod _{p\ll 1}\left (1+\sum _{t\leq k}p^{t(1-1/k)}|G(p^t)|\right )\prod _{p\leq 2Z}\left (1+2k/p\right )
\\ &\ll &\prod _{p\leq 2Z}\left (1+1/p\right )^{2k}
\\ &\ll &(\log Z)^{2k}+1\ll Z^\epsilon 
\end {eqnarray*}
by one of Merten's formulas.  Consequently
\[ \sum _{Z<q\leq 2Z}|G(q)|\ll Z^{1/k-1}\sum _{q\leq 2Z}q^{1-1/k}|G(q)|\ll Z^{1/k-1+\epsilon }.\]
The first ``in particular" claim follows from the main claim after partial summation and then a dyadic split.  The second then 
follows from \eqref {gbound} and the first.  The third follows from the main claim and a dyadic split.  The fourth 
follows from the third and \eqref {gbound}.  The fifth follows from the fourth and partial summation.  The sixth follows from 
\[ \sum _{Z<q\leq 2Z}q^{2-\theta }|G(q)|^2\ll Z^{\epsilon }\]
and a dyadic split, and this in turn follows from the main claim with partial summation.
\end {proof}
\newtheorem {bb}[kfree]{Lemma}
\begin {bb}\label {bb}
Define $H$ as in Definition \ref {qm}.  Define for $q,n\in \mathbb N$
\[ \Phi _q(n)=\sideset {}{'}\sum _{a=1}^q|H(q,a)|^2e\left (\frac {an}{q}\right ),\hspace {5mm}
\Phi _q^*(n)=\sideset {}{'}\sum _{a=1}^q\overline {H(q,a)}e\left (\frac {an}{q}\right )\hspace {5mm}\text {and}\hspace {5mm}
\Phi (q)=\Phi _q(0).\]
(i) Both $\Phi _q(n)$ and $\Phi _q^*(n)$ are, for each $n$, multiplicative in $q$.  If a function $F(q,d)$ defined 
for $q\in \mathbb N$ and $d|q$ satisfies 
for all $(q,q')=1$ and $d|q,d'|q'$ 
\[ F(qq',dd')=F(q,d)F(q',d')\]
then the sum
\[ \sum _{A=1}^qF\Big (q,(q,A)\Big )\Phi _q(A)\]
is multiplcative in $q$.  
\\
\\ (ii) For $q$ a power of a prime and for $d|q$
\[ \sideset {}{'}\sum _{A=1}^{q/d}\Phi _q(-Ad)=\Phi (q)\mu (q/d).\]
(iii) For any $q\in \mathbb N$
\[ \sum _{A=1}^{q}|\Phi _q(A)|\ll q^{1+\epsilon }\]
and the same claim holds with $\Phi _{q}(A)$ replaced by $\Phi _{q}^*(A)$. 
\\
\\ (iv) Let $\eta $ be as in our theorem, $G$ as in Definition \ref {qm}, and define
\[ \rho =\prod _p\left (1-\frac {R_p}{p^k}\right ).\]
Then for any $q\in \mathbb N$ we have
\[ \sum _{a=1}^q\eta (q,a)^2=\frac {\rho ^2}{q}\sum _{d|q}\Phi (d)G(d)^2.\]
\end {bb}
\begin {proof}
For comparison with \cite {kfree} think of $\Phi _q(n)$ as Ramanujan's sum and see Lemma 2.4 of that paper.
\\
\\ (i) All these claims are simple consequences of the fact that $H$ is quasi-multiplcative.
\\
\\ (ii) Write $p$ for the prime in question and suppose $q|p^k$ since otherwise $H(q,a)=0$ so that the claim is trivial.  We have for any $N\in \mathbb N_0$
\begin {equation}\label {mongolia}
\Phi _q(N)=\sum _{n,n'=1}^{R_p}\sideset {}{'}\sum _{a=1}^qe\left (-\frac {a(H_n-H_{n'}-N)}{q}\right )
\end {equation}
so that
\begin {eqnarray*}
\sideset {}{'}\sum _{A=1}^{q/d}\Phi _q(-Ad)&=&\sum _{n,n'=1}^{R_p}
\sideset {}{'}\sum _{a=1}^{q}e\left (-\frac {a(H_n-H_{n'})}{q}\right )
\sideset {}{'}\sum _{A=1}^{q/d}e\left (\frac {-aAd}{q}\right )
\\ &=&\mu (q/d)\sum _{n,n'=1}^{R_p}
\sideset {}{'}\sum _{a=1}^{q}e\left (-\frac {a(H_n-H_{n'})}{q}\right )
\\ &=&\mu (q/d)\Phi _q(0)
\end {eqnarray*}
from \eqref {mongolia}.
\\
\\ (iii) As in (i) the sum in question is 
multiplicative so it is enough to prove the bound for $q$ power of a prime $p$ and as in (ii) it is enough to prove it for $q|p^k$.  In 
that case \eqref {mongolia} implies
\[ \Phi _q(N)=\sum _{n,n'=1}^{R_p}c_q(-H_{n}+H_{n'}+N)\]
so that, since for given $M$ there are only $\ll _r1$ many $(n,n',N)\in \{ 1,...,R_p\} ^2\times \{ 1,...,q\} $ such that $-H_n+H_{n'}+N\equiv M$ modulo $q$,
\[ \sum _{N=1}^q|\Phi _q(N)|\ll \sum _{M=1}^q|c_q(M)|\ll q^{1+\epsilon }\]
by a standard bound for Ramanujan's sum and the proof is similar for $\Phi _{q}^*(n)$.
\\
\\ (iv)  For $P$ a power of a prime $p$ and $A\in \mathbb N$ we have
\begin {equation}\label {HHH}
H(P,A)=\sum _{n=1}^{R_p}e\left (\frac {AH_n}{P}\right )
\end {equation}
where $H_1,...,H_{R_p}$ are the distinct residues represented by $h_1,...,h_r$ modulo $p^k$.  Therefore for $d|q$ with $q$ a power of 
a prime $p$ we have 
\begin {eqnarray}
\sideset {}{'}\sum _{b=1}^{d}H\left (d,b\right )e\left (\frac {-ab}{d}\right )
=
\sum _{n=1}^{R_p}c_{d}(H_n-a).
\end {eqnarray}
By orthogonality and \eqref {jajaja}
\begin {eqnarray*}
g_\mathcal R(q,a)&=&\frac {1}{q}\sum _{b=1}^qG_\mathcal R(q,b)e\left (\frac {ab}{q}\right )
\\ &=&\frac {1}{q}\sum _{d|q}G(d)\sideset {}{'}\sum _{b=1}^{d}H(d,b)e\left (\frac {ab}{d}\right )
\\ &=&\sum _{n=1}^{R_p}\sum _{d|q}G(d)c_d(H_n-a)
\end {eqnarray*}
so that
\begin {eqnarray}\label {wy2}
\sum _{a=1}^q|g_\mathcal S(q,a)|^2&=&\sum _{n,n'=1}^{R_p}\sum _{d,d'|q}G(d)\overline {G(d')}
\sum _{a=1}^qc_d(H_n-a)c_{d'}(a-H_{n'})\notag 
\\ &=&q\sum _{n,n'=1}^{R_p}\sum _{d,d'|q}\frac {G(d)\overline {G(d')}}{[d,d']}
\sum _{a=1}^{[d,d']}c_d(H_n-a)c_{d'}(a-H_{n'}).\hspace {15mm}
\end {eqnarray}
Then
\begin {eqnarray*}
&&\sum _{a=1}^{[d,d']}c_d(H_n-a)c_{d'}(a-H_{n'})
\\ &&\hspace {10mm}=\hspace {4mm}\sideset {}{'}\sum _{A=1}^d\sideset {}{'}\sum _{A'=1}^{d'}e\left (\frac {AH_n}{d}
-\frac {A'H_{n'}}{d'}\right )\sum _{a=1}^{[d,d']}
e\left (\frac {a(-A[d,d']/d+A'[d,d']/d')}{[d,d']}\right )
\\ &&\hspace {10mm}=\hspace {4mm}
[d,d']\sideset {}{'}\sum _{A=1}^d\sideset {}{'}\sum _{A'=1\atop {[d,d']|-A[d,d']/d+A'[d,d']/d'}}^{d'}
e\left (\frac {AH_n}{d}-\frac {A'H_{n'}}{d'}\right )
\end {eqnarray*}
but the only (prime power) $d,d'$ which can 
satisfy these summation conditions are those with $d=d'$, in which case the $A,A'$ sum becomes
\[ \sideset {}{'}\sum _{A=1}^de\left (\frac {A(H_n-H_{n'})}{d}\right )=c_d(H_n-H_{n'})\]
and so
\[ \frac {1}{[d,d']}\sum _{a=1}^{[d,d']}c_d(H_n-a)c_{d'}(H_{n'}-a)=
\left \{ \begin {array}{ll}c_d(H_n-H_{n'})&\text { if }d=d'\\ 0&\text { if not.}\end {array}\right .\]
Therefore \eqref {wy2} says
\begin {eqnarray}\label {lucy}
\sum _{a=1}^q|g_\mathcal R(q,a)|^2&=&q\sum _{n,n'=1}^{R_p}
\sum _{d|q}|G(d)|^2c_d(H_n-H_{n'})\notag 
\\ &=&q\sum _{d|q}|G(d)|^2\sideset {}{'}\sum _{a=1}^q\left |\hspace {1mm}-\sum _{n=1}^{R_p}
e\left (\frac {aH_n}{d}\right )\right |^2\notag 
\\ &=&q\sum _{d|q}\Phi (d)|G(d)|^2.\hspace {15mm}
\end {eqnarray}
from \eqref {HHH}.  This holds initially only for $q$ a prime power,
but the LHS is multiplcative since $g_\mathcal R(q,a)$ is quasi-multiplcative (as in part (i)) and the RHS 
is multiplicative from part (i), so \eqref {lucy} holds in fact for 
general $q$.  Since obviously Lemma \ref {el} says $\eta (q,a)=\rho _\mathcal Rg_\mathcal R(q,a)$ we deduce from \eqref {lucy}
\begin {equation}\label {lucy2}
\sum _{a=1}^q|\eta (q,a)|^2=|\rho _\mathcal R|^2q\sum _{d|q}\Phi (d)|G(d)|^2.
\end {equation}
From \eqref {HHH}  
\[ \sideset {}{'}\sum _{a=1}^q|H(q,a)|^2=\sum _{n,n'=1}^{R_p}c_q(H_n-H_{n'})\]
on prime powers, so that from Definition \eqref {qm}
\begin {equation}\label {drwsnesa}
\sum _{t\geq 0}|G(p^t)|^2\sideset {}{'}\sum _{a=1}^{p^t}|H(p^t,a)|^2=1+
\frac {1/p^{2k}}{(1-R/p^k)^2}\sum _{n,n'=1}^{R_p}\sum _{1\leq t\leq k}c_{p^t}(H_n-H_{n'}).
\end {equation}
But it is easy to establish that for any $D\not |\hspace {1mm}N$
\[ \sum _{d|D}c_d(N)=0\]
and therefore
\[ \sum _{n,n'=1}^{R_p}\sum _{q|p^k\atop {q\not =1}}c_{q}(H_n-H_{n'})=-\sum _{n,n'=1\atop {n\not =n'}}^{R_p}1
+\sum _{n=1}^{R_p}\sum _{q|p^k\atop {q\not =1}}\phi (q)=-R(R-1)+R(p^k-1)=R(p^k-R)\]
which we put in \eqref {drwsnesa} to see that
\[ \sum _{t\geq 0}G(p^t)^2\sideset {}{'}\sum _{a=1}^{p^t}|H(p^t,a)|^2=1+\frac {R(p^k-R)/p^{2k}}{(1-R/p^k)^2}=\frac {1}{1-R/p^k}\]
and so from the multiplicativity of $G$ and from part (i)
\[ \sum _{q=1}^\infty |G(q)|^2\sideset {}{'}\sum _{a=1}^q|H(q,a)|^2=\prod _p\frac {1}{1-R/p^k}=\rho ^{-1}\]
and therefore from \eqref {jajaja}
\[ \sum _{q=1}^\infty \sideset {}{'}\sum _{a=1}^q|G(q,a)|^2=\rho ^{-1}.\]
But from page 92 of \cite {extremal} (more precisely from (1.4), (E) and the following paragraph) the LHS 
is $\rho _\mathcal R^{-1}$ so that in fact $\rho _\mathcal R=\rho $, and the result follows from \eqref {lucy2}.
\end {proof}
%
In the last lemma we showed $\rho _\mathcal R=\rho \not =0$, where $\rho $ is as given in that 
lemma.  From \eqref {exponentialsum2} and \eqref {jajaja} we 
conclude for any $t>0$ and any $q,a\in \mathbb N$ with $(q,a)=1$
\begin {equation}\label {exponentialsum}
\sum _{n\leq t\atop {n\in \mathcal R}}e\left (\frac {an}{q}\right )=\rho G(q)H(q,a)t
+\sum _{\nu =1}^qe\left (\frac {a\nu }{q}\right )E_t(q,\nu )
\end {equation}
which will be our exponential sum approximation in the circle method application. 
\begin {center}
\section {- \hspace {2mm}The circle method application}
\end {center}
In this section we carry out most of the circle method work.
\\
\\ Let $\gamma >0$ be a parameter.  Consider the set of all irreducible fractions in $[0,1]$ with denominator not 
exceeding $\gamma $; the Farey fractions.  If $a'/q'<a/q$ are consecutive Farey fractions in lowest form, define their \emph {median} as 
\[ \frac {a+a'}{q+q'}.\]
Since this lies in $(a'/q',a/q)$ we may partition some\footnote {precisely this interval is $[M,1+M]$ where $X=\lfloor \gamma \rfloor $ 
is the largest permissible denominator and $M$ is the median of $1/X$ and $0$} unit interval $\mathfrak F$ into 
disjoint intervals each containing a Farey point $a/q$ in lowest form and 
extending to the median of $a/q$ with its neighbouring Farey points.  Denoting each interval by $\mathfrak F(a/q)$, the 
\emph {Farey arc at $a/q$}, we see that
\begin {equation}\label {farey999}
\int _\mathfrak Ff(t)dt=\sum _{q\leq \gamma }\sideset {}{'}\sum _{a=1}^q\int _{\mathfrak F(a/q)}f(t)dt
\end {equation}
for any continuous function $f:\mathbb R\rightarrow \mathbb C$.  Denote 
by $\mathfrak U(a/q)$ the interval of unit length centered at $a/q$.  
It can be shown that
\begin {eqnarray}\label {far2}
\left (\frac {a}{q}-\frac {1}{2q\gamma },\frac {a}{q}+\frac {1}{2q\gamma }\right )\subseteq \mathfrak F(a/q)\subseteq 
\left (\frac {a}{q}-\frac {1}{q\gamma },\frac {a}{q}+\frac {1}{q\gamma }\right )\subseteq \mathfrak U(a/q);\hspace {15mm}
\end {eqnarray}
for a discussion of these matters, see Sections 3.1 and 3.8 of \cite {hw}.  
\\ 
\newtheorem {cm}{Lemma}[section]
\begin {cm}\label {cm}
Let $x,\gamma \geq 1$ and $Q>\sqrt x$.  Let $\theta ,\Delta ,G,H$ and $\rho $ be as in Lemma \ref {el}, Definition \ref {qm} and 
Lemma \ref {bb}.  As explained above, denote 
by $\mathfrak F(a/q)$ the Farey arc at $a/q$ in the Farey dissection of order $\gamma $, where $(a,q)=1$, and by $\mathfrak U(a/q)$ the unit 
interval centered at $a/q$.  Define for $t>0$ and $q,a\in \mathbb N$
\[ \Delta _t(q,a)=\sum _{\nu =1}^qe\left (\frac {a\nu }{q}\right )E_t(q,\nu )\]
where $E_t(q,\nu )$ is as in Lemma \ref {el}.  For 
$\alpha \in \mathbb R$ define
\[ f(\alpha )=\sum _{n\leq x\atop {n\in \mathcal R}}e(n\alpha ),\hspace {10mm}
F(\alpha )=\sum _{uv\leq x\atop {u\leq Q}}e(\alpha uv)\hspace {10mm}
\text {and}\hspace {10mm}I(\alpha )=\int _1^xe(\alpha t)dt.\]
For $\alpha \in \mathfrak F(a/q)$ write $\beta =\alpha -a/q$ and define 
\[ J(\alpha )=-2\pi i\beta \int _1^xe(\beta t)\Delta _t(q,a)dt\hspace {10mm}\text {and}\hspace {10mm}
\hat J(\alpha )=e(x\beta )\Delta _x(q,a)+J(\alpha ).\]
For $2\sqrt x\leq \gamma \leq x^{3/4}$ we have
\begin {eqnarray*}
(A)&&\text {for }\alpha \in \mathfrak F(a/q),\hspace {5mm}f(\alpha )=\rho G(q)H(q,a)I(\beta )+\hat J(\alpha )
\\ (B)&&\text {for any }\beta \in \mathbb R\text { and }q\leq x,\hspace {5mm}
\sideset {}{'}\sum _{a=1}^q|H(q,a)|^2F(-a/q-\beta )\ll x
\\ (C)&&\hspace {10mm}\sum _{q\leq 2\sqrt x}\sideset {}{'}\sum _{a=1}^q
\int _{\mathfrak F(a/q)}|F(-\alpha )|\cdot |\hat J(\alpha )|^2d\alpha \ll x^{1+\theta }+\frac {x^{1+2\Delta }}{\gamma }
\\ (D)&&\hspace {10mm}\sum _{q\leq 2\sqrt x}\overline {G(q)}
\sideset {}{'}\sum _{a=1}^q\overline {H(q,a)}\int _{\mathfrak F(a/q)}F(-\alpha )
\overline {I(\beta )}\hat J(\alpha )d\alpha \ll x^{1+\Delta }+\gamma x^{1/2+\theta }
\\ (E)&&\hspace {10mm}\sum _{2\sqrt x<q\leq \gamma }\sideset {}{'}\sum _{a=1}^q
\int _{\mathfrak F(a/q)}|F(-\alpha )|\cdot |f(\alpha )|^2d\alpha \ll x\gamma 
\\ (F)&&\hspace {10mm}
\sideset {}{'}\sum _{a=1}^q\int _{\mathfrak U(a/q)\char92 \mathfrak F(a/q)}
\left |F(-\alpha )\right |\cdot |I(\alpha -a/q)|^2d\alpha \ll 
q^2\gamma ^2.
\end {eqnarray*}
Here the $\ll $ symbol is allowed to contain\footnote {see the notation explained in the introduction} terms of size $x^\epsilon $.
\begin {proof}
This is essentially all contained in \cite {kfree}.  We use a specific notation just for this proof: we 
will write $f(x)\prec \! \! \prec g(x)$ to 
mean $f(x)\ll x^\epsilon g(x)$.  We are basically telling the reader to ignore logs and 
epsilons.  Write $\lambda =1/q\gamma $ so that \eqref {far2} reads
\begin {eqnarray}\label {far}
\left (\frac {a}{q}-\lambda /2,\frac {a}{q}+\lambda /2\right )\subseteq \mathfrak F(a/q)\subseteq 
\left (\frac {a}{q}-\lambda ,\frac {a}{q}+\lambda \right )\subseteq \mathfrak U(a/q).\hspace {15mm}
\end {eqnarray}
For $\alpha \in \mathfrak F(a/q)$ (and so assuming $q\leq \gamma $) we have from \eqref {far} 
that $|\beta |\leq \lambda \leq 1/2\sqrt x$ so from display (2.7), Lemma 2.9, (the second part of) Lemma 2.11 
and Lemma 2.12 of \cite {kfree} we have 
\begin {equation}\label {vaughny}
F(\alpha )\prec \! \! \prec \frac {x}{q\left (1+x|\beta |\right )}+\sqrt x+q
\ll \frac {x}{q}+q,\hspace {4mm}\text { for }\alpha \in \mathfrak F(a/q),\hspace {10mm}
\end {equation}
so that from \eqref {far}
\begin {equation}\label {vaughny3}
\int _{\mathfrak F(a/q)}|F(\alpha )|d\alpha \prec \! \! \prec \frac {1}{q}
\int _{\pm \lambda }\frac {xd\beta }{1+x|\beta |}+\lambda (\sqrt x+q)\prec \! \! \prec \frac {1}{q},
\end {equation}
and 
\begin {equation}\label {vaughny5}
|F(\alpha )|\cdot |\beta |\prec \! \! \prec \frac {1}{q}+\Big (\sqrt x+q\Big )|\beta |\ll \frac {1}{q}
,\hspace {4mm}\text { for }\alpha \in \mathfrak F(a/q).\hspace {10mm}
\end {equation}
From displays (2.7), (2.9), (2.11), Lemma 2.9 and (taking $q=1$ in the first part of) Lemma 2.11 of \cite {kfree}
\begin {equation}\label {vaughny2}
F(\alpha )\ll \sum _{u\leq \sqrt x}\frac {x}{u+x||u\alpha ||}\hspace {10mm}\text { for }\alpha \in \mathbb R. 
\end {equation}
Simply integrating shows
\begin {equation}\label {II}
I(\beta )\ll \frac {x}{1+x|\beta |}.
\end {equation}
Now we prove the claims of the lemma.
\\
\\ (A)  Write $f_t(\alpha )=\Sigma _{n\leq t,n\in \mathcal R}e(n\alpha )$ so that \eqref {exponentialsum} reads 
\begin {eqnarray*}
f_t(a/q)=\rho G(q)H(q,a)t+\Delta _t(q,a)
\end {eqnarray*}
for $t>0$, so that partial summation to
\[ f(\alpha )=\sum _{n\leq x\atop {n\in \mathcal R}}e\left (\frac {an}{q}+\beta n\right )\]
gives
\begin {eqnarray*}
f(\alpha )&=&e(x\beta )f_x(a/q)
-2\pi i\beta \int _1^xe(\beta t)f_t(\alpha )dt
\\ &=&\rho G(q)H(q,a)\left (xe(x\beta )-2\pi i\beta \int _1^xte(\beta t)dt\right )
\\ &&+\hspace  {4mm}e(x\beta )\Delta _x(q,a)-2\pi i\beta \int _1^xe(\beta t)\Delta _t(q,a)dt
\end {eqnarray*}
which gives the result after an integration by parts.  
\\ 
\\ (B)  Take $\Phi _q(n)$ as in Lemma \ref {bb}.  By part (iii) of that lemma
\begin {eqnarray*}
\sideset {}{'}\sum _{a=1}^q|H(q,a)|^2F(-a/q-\beta )&=&\sum _{uv\leq x\atop {u\leq Q}}\Phi _q(-uv)
e\left (-uv\beta \right )
\\ &\ll &x^\epsilon \left (x/q+1\right )\sum _{n=1}^q|\Phi _q(n)|
\\ &\ll &x^\epsilon \left (x/q+1\right )q^{1+\epsilon }.
\end {eqnarray*}
(C)  From \eqref {vaughny5} and \eqref {far} we have for $\alpha \in \mathfrak F(a/q)$
\[ |F(\alpha )|\cdot |J(\alpha )|^2\ll |F(\alpha )|\cdot |\beta |^2\left |\int _1^x\Delta _t(q,a)e(\beta t)dt\right |^2
\prec \! \! \prec \frac {1}{q^2\gamma }\left |\int _1^x\Delta _t(q,a)e(\beta t)dt\right |^2\]
and therefore from \eqref {far}
\begin {eqnarray*}
&&\sideset {}{'}\sum _{a=1}^q\int _{\mathfrak F(a/q)}|F(\alpha )|\cdot |J(\alpha )|^2d\alpha 
\\ &&\hspace {10mm}\prec \! \! \prec \hspace {4mm}\frac {1}{q^2\gamma }\sum _{a=1}^q\int _{\pm \lambda }
\left |\int _1^x\Delta _t(q,a)e(\beta t)dt\right |^2d\beta 
\\ &&\hspace {10mm}=\hspace {4mm}\frac {1}{q^2\gamma }\int _1^x\int _1^x\left (\sum _{a=1}^q\Delta _t(q,a)
\overline {\Delta _{t'}}(q,a)\right )\left (\int _{\pm \lambda }e\Big (\beta (t-t')\Big )d\beta \right )dt'dt.
\end {eqnarray*}
We have ($E_t$ is defined in Lemma \ref {el})
\begin {equation}\label {orthog}
\sum _{a=1}^q\Delta _t(q,a)\overline {\Delta _{t'}}(q,a)=q\sum _{\nu =1}^qE_t(q,\nu )\overline {E_{t'}(q,\nu )}
\end {equation}
and the second factor in the double integral above is
\[ \ll \min \left (\frac {1}{|t-t'|},\lambda \right )\]
so
\begin {eqnarray}\label {can}
&&\sum _{q\leq \gamma }\sideset {}{'}\sum _{a=1}^q\int _{\mathfrak F(a/q)}|F(\alpha )|\cdot |J(\alpha )|^2d\alpha \notag 
\\ &&\hspace {10mm}\prec \! \! \prec \hspace {4mm}
\frac {1}{\gamma }\int _1^x\int _1^x\left (\sum _{q\leq \gamma }\frac {1}{q}\sum _{\nu =1}^q\left |E_t(q,\nu )
\overline {E_{t'}(q,\nu )}\right |\min \left (\frac {1}{|t-t'|},1\right )\right )dt'dt\notag 
\\ &&\hspace {10mm}=:\hspace {4mm}\frac {\mathcal V(x,\gamma )}{\gamma }.
\end {eqnarray}
Applying twice the Cauchy-Schwarz inequality we see that
\begin {eqnarray}\label {ilka}
\mathcal V(x,\gamma )&\leq &\int _1^x\left (\sum _{q\leq \gamma }\frac {1}{q}\sum _{\nu =1}^q
\left |E_t(q,\nu )\right |^2\int _1^x\min \left (\frac {1}{|t-t'|},1\right )dt'\right )dt \notag 
\\ &\prec \! \! \prec &x\cdot \max _{1\leq t\leq x}\left (\sum _{q\leq \gamma }\frac {1}{q}\sum _{\nu =1}^q\left |E_t(q,\nu )\right |^2\right )
\end {eqnarray}
therefore from \eqref {can}
\begin {eqnarray}\label {can2}
\sum _{q\leq \gamma }\sideset {}{'}\sum _{a=1}^q\int _{\mathfrak F(a/q)}|F(\alpha )|\cdot |J(\alpha )|^2d\alpha 
\prec \! \! \prec \frac {x}{\gamma }\cdot \max _{1\leq t\leq x}\left (\sum _{q\leq \gamma }
\frac {1}{q}\sum _{\nu =1}^q\left |E_t(q,\nu )\right |^2\right ).
\end {eqnarray}
From \eqref {vaughny3} and \eqref {orthog}
\[ \sideset {}{'}\sum _{a=1}^q|\Delta _x(q,a)|^2\int _{\mathfrak F(a/q)}|F(\alpha )|d\alpha 
\prec \! \! \prec \sum _{\nu =1}^q|E_x(q,\nu )|^2\]
therefore
\begin {eqnarray*}
\sum _{q\leq 2\sqrt x}\sideset {}{'}\sum _{a=1}^q|\Delta _x(q,a)|^2\int _{\mathfrak F(a/q)}|F(\alpha )|d\alpha  
\prec \! \! \prec \max _{1\leq t\leq x}\left (\sum _{q\leq 2\sqrt x}\sum _{\nu =1}^q|E_t(q,\nu )|^2\right )
\end {eqnarray*}
which with \eqref {can2} says
\begin {eqnarray*}
&&\sum _{q\leq 2\sqrt x}\sideset {}{'}\sum _{a=1}^q\int _{\mathfrak F(a/q)}|F(\alpha )|\cdot |\hat J(\alpha )|^2d\alpha 
\\ &&\hspace {15mm}\prec \! \! \prec \hspace {4mm}\max _{1\leq t\leq x}\left (\sum _{q\leq 2\sqrt x}
\sum _{\nu =1}^q|E_t(q,\nu )|^2+\frac {x}{\gamma }\sum _{q\leq \gamma }\frac {1}{q}\sum _{\nu =1}^q|E_t(q,\nu )|^2\right )
\\ &&\hspace {15mm}\prec \! \! \prec \hspace {4mm}
x^{1+\theta }+x^{2\Delta }+x^{1/2+\Delta }
+\frac {x}{\gamma }\Big (\gamma ^{1-2\theta }x^{2\theta }+x^{2\Delta }\Big )
\end {eqnarray*}
from Lemma \ref {el} (ii).  The second term is bounded by the fifth (since $\gamma \leq x$), and 
the third and fourth are bounded by the first (the third since $\Delta \leq 1/2$ unless $k=2$, in 
which case $1/2+\Delta \leq 3/2=1+\theta $, and the fourth since $\gamma \geq \sqrt x$).
\\ 
\\ (D)  For $\alpha \in \mathfrak F(a/q)$ write
\[ C(\alpha )=G(q)H(q,a)I(\beta ).\]
For $\alpha \in \mathfrak F(a/q)\char92 (a/q-\lambda /2,a/q+\lambda /2)$ we have $|\beta |\gg \lambda $ so 
from \eqref {vaughny5} and \eqref {II}
\[ F(-\alpha )|I(\beta )|^2\prec \! \! \prec \frac {1}{q|\beta |}\left (\frac {x}{1+x|\beta |}\right )^2
\ll \frac {1}{q\lambda ^3}\]
so that from \eqref {far}
\[ \int _{\mathfrak F(a/q)\char92 (a/q-\lambda /2,a/q+\lambda /2)}|F(-\alpha )|\cdot |I(\beta )|^2d\alpha 
\prec \! \! \prec \frac {1}{q\lambda ^2}=q\gamma ^2\]
and therefore by \eqref {hbound}
\begin {eqnarray*}
\sum _{q\leq 2\sqrt x}\sideset {}{'}\sum _{a=1}^q\int _{\mathfrak F(a/q)\char92 (a/q-\lambda /2,a/q+\lambda /2)}
|F(-\alpha )|\cdot |C(\alpha )|^2d\alpha &\prec \! \! \prec &\gamma ^2\sum _{q\leq 2\sqrt x}q^2|G(q)|^2
\\ &\prec \! \! \prec &\gamma ^2x^{\theta /2}
\end {eqnarray*}
from Lemma \ref {cc}.  Therefore the Cauchy-Schwarz Inequality and part (C) imply
\begin {eqnarray*}
&&\left (\sum _{q\leq 2\sqrt x}\sideset {}{'}\sum _{a=1}^q\int _{\mathfrak F(a/q)\char92 (a/q-\lambda /2,a/q+\lambda /2)}
|F(-\alpha )|\cdot |C(\alpha )|\cdot |\hat J(\alpha )|d\alpha \right )^2
\\ &&\hspace {10mm}\leq \hspace {4mm}\left (\sum _{q\leq 2\sqrt x}\sideset {}{'}\sum _{a=1}^q
\int _{\mathfrak F(a/q)\char92 (a/q-\lambda /2,a/q+\lambda /2)}
|F(-\alpha )|\cdot |C(\alpha )|^2d\alpha \right )
\\ &&\hspace {14mm}\times \hspace {4mm}
\left (\sum _{q\leq 2\sqrt x}\sideset {}{'}\sum _{a=1}^q\int _{\mathfrak F(a/q)}
|F(-\alpha )|\cdot |\hat J(\alpha )|^2d\alpha \right )
\\ &&\hspace {10mm}\prec \! \! \prec \hspace {4mm}
\gamma ^2x^{\theta /2}\left (x^{1+\theta }+\frac {x^{1+2\Delta }}{\gamma }\right )
\end {eqnarray*}
so that, since $\gamma \leq x^{3/4}$,
\begin {eqnarray}\label {banana}
&&\sum _{q\leq 2\sqrt x}\sideset {}{'}\sum _{a=1}^q\int _{\mathfrak F(a/q)\char92 (a/q-\lambda /2,a/q+\lambda /2)}
|F(-\alpha )|\cdot |C(\alpha )|\cdot |\hat J(\alpha )|d\alpha \notag 
\\ &&\hspace {10mm}\prec \! \! \prec \hspace {4mm}\gamma x^{1/2+3\theta /4}+\gamma ^{1/2}x^{1/2+\Delta +\theta /4}\notag 
\\ &&\hspace {10mm}\ll \hspace {4mm}\gamma x^{1/2+\theta }+x^{1+\Delta }.
\end {eqnarray}
For $\alpha \in \mathfrak F(a/q)$ (so assuming $(a,q)=1$ and $q\leq \gamma $) and $u\leq \sqrt x$ we have from \eqref {far}
\[ |u\beta |\leq \frac {1}{2q}\]
so if $q\not |\hspace {1mm}u$ then
\[ ||u\alpha ||\geq \left |\left |\frac {ua}{q}\right |\right |-||u\beta ||\gg \left |\left |\frac {ua}{q}\right |\right |\]
and therefore using the standard bound for a linear exponential sum
\begin {eqnarray}\label {adio2}
H_q(\alpha )&:=&\sum _{u\leq \sqrt x\atop {q\not |\hspace {1mm}u}}\left (\sum _{v\leq x/u}
+\sum _{\sqrt x<u\leq Q,x/v}\right )e(\alpha uv)\notag 
\\ &\ll &\sum _{u\leq \sqrt x\atop {q\not |\hspace {1mm}u}}\frac {1}{||ua/q||}\notag 
\\ &\ll &(\sqrt x/q+1)\sum _{u=1}^{q-1}\frac {1}{||ua/q||}\notag 
\\ &\ll &\log q\Big (\sqrt x+q\Big )\ll \gamma \log q
\end {eqnarray}
so that, breaking the $u$ summation in the definition of $F$ at $\sqrt x$ and then swapping sums in the second part,
\begin {eqnarray}\label {adio}
F(\alpha )&=&\sum _{u\leq \sqrt x}\sum _{v\leq x/u}e(\alpha uv)
+\sum _{v\leq \sqrt x}\sum _{\sqrt x<u\leq Q,x/v}e(\alpha uv)\notag 
\\ &=&\sum _{u\leq \sqrt x\atop {q|u}}\sum _{v\leq x/u}e(\beta uv)
+\sum _{u\leq \sqrt x\atop {q|u}}\sum _{\sqrt x<v\leq Q,x/u}e(\beta uv)+H_q(\alpha )\notag 
\\ &=:&K_q(\beta )+\mathcal O\left (\gamma \log q\right )
\end {eqnarray}
whenever $\alpha \in \mathfrak F(a/q)$.  Therefore by \eqref {far} and \eqref {hbound}
\begin {eqnarray*}
&&\sideset {}{'}\sum _{a=1}^q\overline {H(q,a)}\int _{a/q\pm \lambda /2}F(-\alpha )\overline {I(\beta )}\hat J(\alpha )d\alpha 
\\ &&\hspace {10mm}=\hspace {4mm}\int _{\pm \lambda /2}
K_q(-\beta )I(-\beta )\left (\sideset {}{'}\sum _{a=1}^q\overline {H(q,a)}
\hat J(a/q+\beta )\right )d\beta 
+\mathcal O\left (q^\epsilon \gamma \sideset {}{'}\sum _{a=1}^q\int _{\mathfrak F(a/q)}
\left |I(\beta )\hat J(\alpha )\right |d\alpha \right )
\\ &&\hspace {10mm}=:\hspace {4mm}\int _{\pm \lambda /2}K_q(-\beta )A_q(\beta )d\beta 
+\mathcal O\left (q^\epsilon \gamma \sideset {}{'}\sum _{a=1}^q\int _{\mathfrak F(a/q)}|B(\alpha )|d\alpha \right )
\end {eqnarray*}
so that
\begin {eqnarray*}
&&\sum _{q\leq 2\sqrt x}\sideset {}{'}\sum _{a=1}^q\int _{a/q\pm \lambda /2}F(-\alpha )\overline {C(\alpha )}\hat J(\alpha )d\alpha 
\\ &&\hspace {10mm}\prec \! \! \prec \hspace {4mm}
\sum _{q\leq 2\sqrt x}|G(q)|\left |\int _{\pm \lambda /2}K_q(-\beta )A_q(\beta )d\beta \right |
+\gamma \sum _{q\leq 2\sqrt x}|G(q)|\sideset {}{'}\sum _{a=1}^q\int _{\mathfrak F(a/q)}|B(\alpha )|d\alpha 
\end {eqnarray*}
and therefore from \eqref {banana}
\begin {eqnarray}\label {everyday}
&&\sum _{q\leq 2\sqrt x}\overline {G(q)}\sideset {}{'}\sum _{a=1}^q\overline {H(q,a)}\int _{\mathfrak F(a/q)}
F(-\alpha )\overline {I(\beta )}\hat J(\alpha )d\alpha \notag 
\\ &&\hspace {10mm}=\hspace {4mm}\sum _{q\leq 2\sqrt x}\sideset {}{'}\sum _{a=1}^q
\int _{a/q\pm \lambda /2}F(-\alpha )\overline {C(\alpha )}\hat J(\alpha )d\alpha \notag  
\\ &&\hspace {14mm}+\hspace {4mm}\mathcal O\left (\sum _{q\leq 2\sqrt x}\sideset {}{'}
\sum _{a=1}^q\int _{\mathfrak F(a/q)\char92 (a/q-\lambda /2,a/q+\lambda /2)}
|F(-\alpha )|\cdot |C(\alpha )|\cdot |\hat J(\alpha )|d\alpha \right )\notag 
\\ &&\hspace {10mm}\prec \! \! \prec \hspace {4mm}
\sum _{q\leq \gamma }|G(q)|\left |\int _{\pm \lambda /2}K_q(-\beta )A_q(\beta )d\beta \right |\notag 
\\ &&\hspace {14mm}+\hspace {4mm}\gamma \sum _{q\leq \gamma }|G(q)|\sideset {}{'}\sum _{a=1}^q
\int _{\mathfrak F(a/q)}|B(\alpha )|d\alpha +\gamma x^{1/2+\theta }+x^{1+\Delta }.\hspace {15mm}
\end {eqnarray}
Recall the definition of $E_t(q,\nu )$ from Lemma \ref {el} and of $H(q,a)$ from 
Lemma \ref {bb}, and let $\Phi ^*(n)$ be as in Lemma \ref {bb}.  From Lemma \ref {el} (i) and 
then Lemma \ref {bb} (iii) we have for $q,t\leq x$
\begin {eqnarray*}
\sideset {}{'}\sum _{a=1}^q\Delta _t(q,a)\overline {H(q,a)}&=&\sum _{\nu =1}^q\Phi _q^*(\nu )\left (
\sum _{n\leq t\atop {n\in \mathcal R\atop {n\equiv \nu (q)}}}1
-t\eta (q,v)\right )
\\ &\prec \! \! \prec &t^\Delta \sum _{\nu =1}^q|\Phi _q^*(\nu )|
\\ &\prec \! \! \prec &qt^\Delta 
\end {eqnarray*}
so that for any $\beta \in \mathbb R$
\begin {eqnarray*}
\sideset {}{'}\sum _{a=1}^q\overline {H(q,a)}\hat J(a/q+\beta )&\ll &
\left |\sideset {}{'}\sum _{a=1}^q\Delta _x(q,a)\overline {H(q,a)}\right |
+|\beta |\int _1^x\left |\sideset {}{'}\sum _{a=1}^q\Delta _t(q,a)\overline {H(q,a)}\right |dt\notag 
\\ &\prec \! \! \prec &qx^{\Delta }(1+|\beta |x)
\end {eqnarray*}
and therefore from \eqref {II}
\begin {eqnarray}\label {llaeth}
A_q(\beta )\prec \! \! \prec qx^{1+\Delta }.
\end {eqnarray}
Therefore from \eqref {adio}, \eqref {far} and \eqref {vaughny3}
\[ \int _{\pm \lambda /2}K_q(-\beta )A_q(\beta )d\beta \prec \! \! \prec qx^{1+\Delta }
\int _{\pm 1/2q\gamma }\Big (|F(-a/q-\beta )|+\gamma \Big )d\beta \prec \! \! \prec 
qx^{1+\Delta }\left (\frac {1}{q}+\gamma \lambda \right )
\ll x^{1+\Delta }\]
and so from \eqref {gbound}
\begin {eqnarray}\label {lovesong}
\sum _{q\leq \gamma }|G(q)|\left |\int _{\pm \lambda /2}K_q(-\beta )A_q(\beta )d\beta \right |\prec \! \! \prec x^{1+\Delta }.
\end {eqnarray}
We have from \eqref {far}
\begin {eqnarray*}
&&\sideset {}{'}\sum _{a=1}^q\int _{\mathfrak F(a/q)}|J(\alpha )|^2d\alpha 
\\ &&\hspace {10mm}\leq \hspace {4mm}\int _1^x\int _1^x\left (\sum _{a=1}^q\Delta _t(q,a)\overline {\Delta _{t'}(q,a)}
\right )\cdot \left |\int _{\pm \lambda }|\beta |^2e\Big (\beta (t-t')\Big )d\beta \right |dt'dt.
\end {eqnarray*}
The first factor is from \eqref {orthog}
\[ q\sum _{\nu =1}^qE_t(q,\nu )\overline {E_{t'}(q,\nu )}\]
and the second factor is
\[ \ll \min \left (\frac {\lambda ^2}{t-t'},\lambda ^3\right )\]
so 
\begin {eqnarray}\label {ilka3}
&&\sum _{q\leq \gamma }\frac {1}{q^{1-\theta }}\sideset {}{'}\sum _{a=1}^q\int _{\mathfrak F(a/q)}|J(\alpha )|^2d\alpha \notag 
\\ &&\hspace {10mm}\ll \hspace {4mm}\int _1^x\int _1^x\left (\sum _{q\leq \gamma }q^{\theta }\lambda ^2
\sum _{\nu =1}^q\left |E_t(q,\nu )\overline {E_{t'}(q,\nu )}\right |
\min \left (\frac {1}{t-t'},1\right )\right )dt'dt\notag 
\\ &&\hspace {10mm}=:\hspace {4mm}\mathcal U(x,\gamma ).
\end {eqnarray}
As in \eqref {ilka} we have
\begin {eqnarray*}
\mathcal U(x,\gamma )&\leq &\int _1^x\left (\sum _{q\leq \gamma }q^\theta \lambda ^2\sum _{\nu =1}^q
\left |E_t(q,\nu )\right |^2\int _1^x\min \left (\frac {1}{t-t'},1\right )dt'\right )dt
\\ &\prec \! \! \prec &\frac {x}{\gamma ^2}\cdot \max _{1\leq t\leq x}\left (
\sum _{q\leq \gamma }q^{\theta -2}\sum _{\nu =1}^q\left |E_t(q,\nu )\right |^2\right )
\end {eqnarray*}
so that \eqref {ilka3} says
\begin {eqnarray}\label {ilka66}
\sum _{q\leq \gamma }\frac {1}{q^{1-\theta }}\sideset {}{'}\sum _{a=1}^q\int _{\mathfrak F(a/q)}|J(\alpha )|^2d\alpha 
\prec \! \! \prec \frac {x}{\gamma ^2}\cdot \max _{1\leq t\leq x}\left (
\sum _{q\leq \gamma }\frac {1}{q^{2-\theta }}\sum _{\nu =1}^q\left |E_t(q,\nu )\right |^2\right ).\hspace {10mm}
\end {eqnarray}
From \eqref {far} and \eqref {orthog} 
\[ \sideset {}{'}\sum _{a=1}^q|\Delta _x(a/q)|^2\int _{\mathfrak F(a/q)}d\alpha 
\prec \! \! \prec q\lambda \sum _{\nu =1}^q|E_x(q,\nu )|^2=\frac {1}{\gamma }\sum _{\nu =1}^q|E_x(q,\nu )|^2\]
so that
\begin {eqnarray*}
\sum _{q\leq \gamma }\frac {1}{q^{1-\theta }}\sideset {}{'}\sum _{a=1}^q\int _{\mathfrak F(a/q)}|\Delta (a/q)|^2d\alpha 
\prec \! \! \prec \frac {1}{\gamma }\cdot \max _{1\leq t\leq x}
\left (\sum _{q\leq \gamma }\frac {1}{q^{1-\theta }}\sum _{\nu =1}^q|E_t(q,\nu )|^2\right )
\end {eqnarray*}
and therefore from \eqref {ilka66} and Lemma \ref {el} (ii)
\begin {eqnarray}\label {georgia1}
&&\sum _{q\leq \gamma }\frac {1}{q^{1-\theta }}\sideset {}{'}\sum _{a=1}^q\int _{\mathfrak F(a/q)}|\hat J(\alpha )|^2d\alpha \notag 
\\ &&\hspace {10mm}\prec \! \! \prec \hspace {4mm}\frac {1}{\gamma }\cdot \max _{1\leq t\leq x}
\left (\sum _{q\leq \gamma }\frac {1}{q^{1-\theta }}\sum _{\nu =1}^q|E_t(q,\nu )|^2
+\frac {x}{\gamma }\sum _{q\leq \gamma }\frac {1}{q^{2-\theta }}\sum _{\nu =1}^q\left |E_t(q,\nu )\right |^2\right )\notag 
\\ &&\hspace {10mm}\prec \! \! \prec \hspace {4mm}\frac {1}{\gamma }\left (\gamma ^{1-\theta }x^{2\theta }+x^{2\Delta }
+\gamma ^\theta x^\Delta +\frac {x^{1+2\Delta }}{\gamma }\right )\notag 
\\ &&\hspace {10mm}\ll \hspace {4mm}x^{2\theta }+\frac {x^{1+2\Delta }}{\gamma ^2},
\end {eqnarray}
the second term being less than the fourth in the punultimate line (since $x\geq \gamma $), and 
the third less than the first (since $2\theta \geq \Delta $ and $1-\theta \geq \theta $).  From orthogonality
\[ \int _{\pm \lambda }|I(\beta )|^2\ll x\]
so from Lemma \ref {cc}
\begin {equation}\label {georgia11}
\sum _{q\leq \gamma }q^{1-\theta }|G(q)|^2
\sum _{a=1}^q\int _{\pm \lambda }|I(\beta )|^2d\beta \prec \! \! \prec x.
\end {equation}
From the Cauchy-Schwarz inequality and then \eqref {georgia1} and \eqref {georgia11}
\begin {eqnarray}\label {robin}
&&\sum _{q\leq \gamma }|G(q)|\sideset {}{'}\sum _{a=1}^q\int _{\mathfrak F(a/q)}|B_q(\alpha )|d\alpha \notag 
\\ &&\hspace {10mm}\leq \hspace {4mm}\left (\sum _{q\leq \gamma }q^{1-\theta }|G(q)|^2
\sum _{a=1}^q\int _{\pm \lambda }|I(\beta )|^2d\beta \right )^{1/2}
\left (\sum _{q\leq \gamma }\frac {1}{q^{1-\theta }}
\sideset {}{'}\sum _{a=1}^q\int _{\mathfrak F(a/q)}|\hat J(\alpha )|^2d\alpha \right )^{1/2}\notag 
\\ &&\hspace {10mm}\prec \! \! \prec \hspace {4mm}x^{1/2+\theta }+\frac {x^{1+\Delta }}{\gamma }.
\end {eqnarray}
From this, \eqref {everyday} and \eqref {lovesong} we deduce 
\begin {eqnarray*}
&&\sum _{q\leq 2\sqrt x}\overline {G(q)}\sideset {}{'}\sum _{a=1}^q\overline {H(q,a)}\int _{\mathfrak F(a/q)}
F(-\alpha )\overline {I(\beta )}\hat J(\alpha )d\alpha 
\prec \! \! \prec \gamma \left (x^{1/2+\theta }+\frac {x^{1+\Delta }}{\gamma }\right ).
\end {eqnarray*}
(E)  
For $q\gg \sqrt x$ and $\alpha \in \mathfrak F(a/q)$ (so assuming $q\leq \gamma )$ we have from \eqref {vaughny}
\[ F(\alpha )\ll \gamma \]
therefore (since the collection of all Farey arcs gives some interval of unit length)
\begin {eqnarray}\label {Fminor}
\sum _{2\sqrt x<q\leq \gamma }\sideset {}{'}\sum _{a=1}^q\int _{\mathfrak F(a/q)}|F(-\alpha )|
\cdot |f(\alpha )|^2d\alpha \ll \gamma \int _0^1|f(\alpha )|^2\ll x\gamma 
\end {eqnarray}
by orthogonality.
\\
\\ (F)  Write
\[ \mathcal V=\sideset {}{'}\sum _{a=1}^q\sum _{u\leq \sqrt x}\int _{u/2q\gamma <|\beta |\leq 1/2}
\frac {xd\beta }{|\beta |^2\left (u+x||ua/q+u\beta ||\right )}\]
so that \eqref {vaughny2} and the bound $I(\beta )\ll 1/|\beta |$ (from \eqref {II}) says
\begin {equation}\label {vodan69}
\sideset {}{'}\sum _{a=1}^q\int _{1/2q\gamma <|\beta |\leq 1/2}|F(-a/q-\beta )|\cdot |I(\beta )|^2d\beta \prec \! \! \prec \mathcal V.
\end {equation}
We have
\begin {eqnarray*}\int _{1/2q\gamma <|\beta |\leq 1/2}\frac {xd\beta }{\beta ^2\left (u+x||ua/q+u\beta ||\right )}
&=&u\int _{u/2q\gamma <|t|\leq u/2}\frac {xdt}{t^2\left (u+x||ua/q+t||\right )}
\end {eqnarray*}
and the part of the integral with $t\geq 1/2$ is
\begin {eqnarray*}
&\leq &\sum _{0<|j|\leq u/2}\int _{-1/2}^{1/2}\frac {xdt}{(j+t)^2\left (u+x||ua/q+j+t||\right )}
\\ &\ll &\left (\int _{|ua/q+t|\leq 1/x}+\int _{1/x\leq |ua/q+t|\leq 1/2}\right )\frac {xdt}{u+x||ua/q+t||}
\\ &\ll &1+\log x\ll x^\epsilon 
\end {eqnarray*}
so that the whole integral in the definition of $\mathcal V$ is
\[ u\int _{u/2q\gamma <|t|\leq 1/2}\frac {xdt}{t^2\left (u+x||ua/q+t||\right )}+\mathcal O(x^\epsilon u)\]
and therefore
\begin {eqnarray*}
\mathcal V&=&\sum _{u\leq \sqrt x}u\int _{u/2q\gamma <|t|\leq 1/2}\frac {1}{t^2}
\left (\sideset {}{'}\sum _{a=1}^q\frac {x}{u+x||ua/q+t||}\right )dt
+\mathcal O\left (x^\epsilon \sideset {}{'}\sum _{a=1}^q\sum _{u\leq \sqrt x}u\right )
\\ &=&\sum _{u\leq \sqrt x}u\int _{u/2q\gamma <|t|\leq 1/2}\frac {1}{t^2}
\left (\sideset {}{'}\sum _{a=1}^q\frac {x}{u+x||ua/q+t||}\right )dt
+\mathcal O\left (x^{1+\epsilon }q\right ).
\end {eqnarray*}
Write $q'=q/(q,u)$ and $u'=u/(q,u)$.  The inner sum is
\begin {eqnarray*}
&\leq &(q,u)\sum _{a=1}^{q'}\frac {x}{u+x||u'a/q'+t||}
\\ &=&(q,u)\sum _{a=1}^{q'}\frac {x}{u+x||a/q'+t||}
\\ &\ll &(q,u)\sum _{a=1\atop {||a/q'+t||\leq 1/2q'}}^{q'}\frac {x}{u+x||a/q'+t||}+(q,u)q'\log q'
\end {eqnarray*}
so that
\begin {eqnarray*}
\mathcal V&\prec \! \! \prec &\sum _{u\leq \sqrt x}u(q,u)\int _{u/2q\gamma <|t|\leq 1/2}\frac {1}{t^2}
\left (\sum _{a=1\atop {||a/q'+t||\leq 1/2q'}}^{q'}\frac {x}{u+x||a/q'+t||}\right )dt
\\ &&\hspace {10mm}+\hspace {4mm}
q\sum _{u\leq \sqrt x}u\int _{u/q\gamma <|t|\leq 1/2}\frac {dt}{t^2}\hspace {2mm}+\hspace {2mm}xq
\\ &\ll &\sum _{u\leq \sqrt x}u(q,u)\int _{u/2q\gamma <|t|\leq 1/2}
\left (\sum _{a=1\atop {||a/q'+t||\leq 1/2q'}}^{q'}\frac {x}{t^2(u+x||a/q'+t||)}\right )dt
+q^2\gamma \sqrt x
\end {eqnarray*}
In the sum we have $t\gg a/q'$ so that the whole integral is for $q\leq x$
\begin {eqnarray*}
&\ll &{q'}^2\int _{u/2q\gamma <|t|\leq 1/2}\left (\sum _{a=1}^{q'}\frac {x}{a^2(u+x||a/q'+t||)}\right )dt\prec \! \! \prec {q'}^2
\end {eqnarray*}
and therefore
\begin {eqnarray*}
\mathcal V\prec \! \! \prec \sum _{u\leq \sqrt x}u(q,u){q'}^2+q^2\gamma \sqrt x\ll q^2\gamma \sqrt x
\end {eqnarray*}
so that from \eqref {vodan69} 
\[ \sideset {}{'}\sum _{a=1}^q\int _{\lambda /2<|\beta |\leq 1/2}|F(-a/q-\beta )|\cdot |I(\beta )|^2d\beta 
\prec \! \! \prec q^2\gamma ^2\]
or in other words
\[ \sideset {}{'}\sum _{a=1}^q\int _{(-1/2,1/2)\char92 X}|F(-a/q-\beta )|\cdot |I(\beta )|^2d\alpha \prec \! \! \prec q^2\gamma ^2\]
for any subset 
\[ \left (-\frac {\lambda }{2},\frac {\lambda }{2}\right )
\subseteq X\subseteq (-1/2,1/2).\]
The result now follows from \eqref {far}.
\end {proof}
\end {cm}
\begin {center}
\section {- \hspace {2mm}Evaluation of a character sum}
\end {center}
In the last stage of the proof we will be left with a quantity which we have chosen to analyse with Perron's formula.  The main 
difficulty will be evaluating 
\[ \sum _{n\leq X}\frac {\chi (n)}{n^s}\]
where $\chi $ is a Dirichlet character and $s=it$ for $t\in \mathbb R$.
\\
\\ As in Chapter 9 of \cite {montgomeryvaughan} we make the convention that a primitive character may be 
principal (and so necessarily of modulus one).
\\
\newtheorem {mvt}{Lemma}[section]
\begin {mvt}\label {mvt}
For any $M\in \mathbb N$, $Q>0$, $t_0\in \mathbb R$, $T\geq 1$, and any primitive character $\chi $ modulo $M$, 
\[ \int _1^T\frac {L(iv+it_0,\chi )Q^{iv}dv}{v}\ll \sqrt M\Big (1+\sqrt {|t_0|}\Big ).\]
Here the $\ll $ may contain\footnote {as explained in the introduction} terms up to $M^\epsilon ,T^\epsilon ,|t_0|^\epsilon $.
\end {mvt}
\begin {proof}
Throughout we allow the $M^\epsilon ,T^\epsilon ,|t_0|^\epsilon $ terms to go into the $\ll ,\mathcal O$ symbols - we are basically 
telling the reader to ignore logs and epsilons.
\\
\\ We first suppose $\chi $ is non-principal.  Take parameters $Z_0>Z\geq 2$.  Summing 
by parts and applying the Polya-Vinogradov 
Inequality (Theorem 9.18 of \cite {montgomeryvaughan}) we have for any $t\in \mathbb R$
\begin {eqnarray*}
\sum _{Z<n\leq Z_0}\frac {\overline \chi (n)}{n^{1-it}}&=&\frac {1}{Z^{1-it}}\sum _{Z<n\leq Z_0}\overline \chi (n)
+(1-it)\int _Z^{Z_0}\frac {1}{y^{2-it}}\left (\sum _{Z<n\leq y}\overline \chi (n)\right )dy
\\ &\ll &\frac {(1+|t|)\sqrt M}{Z}
\end {eqnarray*}
so that letting $Z_0\rightarrow \infty $ 
\begin {equation}\label {pnawn}
L(1-it,\overline \chi )=\sum _{n\leq Z}\frac {\overline \chi (n)}{n^{1-it}}
+\mathcal O\left (\frac {\log Z}{1+|t|}\right )
\end {equation}
so long as 
\begin {equation}\label {goro1}
Z>(1+|t|)^2\sqrt M;
\end {equation}
on the other hand Theorem 4.11 of \cite {titchmarsh} says that if $Z>1+|t|$ then
\[ \zeta (1-it)=\sum _{n\leq Z}\frac {1}{n^{1-it}}+\mathcal O\left (\frac {\log Z}{1+|t|}\right )\]
so that \eqref {pnawn} subject to \eqref {goro1} remains true also in the case of principal $\chi $ (that is, $M=1$).  For 
any $\kappa $ there is some $A_\kappa $ for which
\[ \sin \left (\frac {\pi (it+\kappa )}{2}\right )=
A_\kappa e^{\pi |t|/2}\left (1+\mathcal O\left (\frac {1}{1+|t|}\right )\right )\]
and by standard formulas for the Gamma function there is some $B$ for which
\[ \Gamma (1-it)=B|t|^{1/2-it}e^{-\pi |t|/2+it}\left (1+\mathcal O\left (\frac {1}{1+|t|}\right )\right )\]
so that with \eqref {pnawn} we have for any $t\in \mathbb R$
\begin {eqnarray*}
&&L(1-it,\overline \chi )\Gamma (1-it)\sin \left (\frac {\pi (it+\kappa )}{2}\right )
\\ &&\hspace {15mm}=\hspace {4mm}\left (\sum _{n\leq Z}\frac {\overline \chi (n)}{n^{1-it}}+\mathcal O\left (\frac {\log Z}{1+|t|}\right )\right )
\left (BA_\kappa |t|^{1/2-it}e^{it}+\mathcal O\left (\frac {1}{(1+|t|)^{1/2}}\right )\right )
\\ &&\hspace {15mm}=\hspace {4mm}BA_\kappa |t|^{1/2-it}e^{it}\sum _{n\leq Z}\frac {\overline \chi (n)}{n^{1-it}}
+\mathcal O\left (\log Z\right )
\end {eqnarray*}
so long as \eqref {goro1} holds.  Let $\kappa $ and $\epsilon (\chi )$ be given respectively as 
in (10.15) and (10.17) of \cite {montgomeryvaughan}; from the comments immediately following (10.17) we 
have $\epsilon (\chi )\ll 1 $.  Therefore Corollary 10.9 of \cite {montgomeryvaughan} and the last equality say that for some $C_{\chi }\ll 1$ we have 
for any $t\in \mathbb R$
\begin {eqnarray*}
L(it,\chi )&=&\pi ^{-1}\epsilon (\chi )\sqrt M\left (\frac {2\pi }{M}\right )^{it}
L(1-it,\overline \chi )\Gamma (1-it)\sin \left (\frac {\pi (it+\kappa )}{2}\right )
\\ &=&C_\chi \sqrt M\left (\frac {2\pi }{M}\right )^{it}|t|^{1/2-it}e^{it}\sum _{n\leq Z}\frac {\overline \chi (n)}{n^{1-it}}
+\mathcal O\left (\sqrt M\log Z\right )
\\ &=&C_\chi \sqrt {M|t|}\sum _{n\leq Z}\frac {\overline \chi (n)e^{it(\log (2\pi /M)-\log |t|+1+\log n)}}{n}
+\mathcal O\left (\sqrt M\log Z\right )
\end {eqnarray*}
so long as \eqref {goro1} holds.  Therefore for any $1\leq v\leq T$ and so long as 
\begin {equation}\label {goro}
Z>(1+T+|t_0|)^2\sqrt M
\end {equation}
we have
\begin {eqnarray*}
\frac {L\left (i(v+t_0),\chi \right )Q^{iv}}{v}&=&\frac {C_\chi \sqrt {M|v+t_0|}}{v}
\sum _{n\leq Z}\frac {\overline \chi (n)e\left (f(v)\right )}{n}
+\mathcal O\left (\frac {\sqrt M\log Z}{v}\right ),
\end {eqnarray*}
where
\[ f(v)=f_{n,M,Q,t_0}(v)=\frac {(v+t_0)(\log (2\pi /M)-\log |v+t_0|+1+\log n)+v\log Q}{2\pi };\]
note that $f$ is twice differentiable for $v+t_0\not =0$ and there we have
\begin {equation}\label {deriv}
f''(v)=\pm \frac {1}{2\pi |v+t_0|}. 
\end {equation}
Therefore
\begin {eqnarray}\label {bora}
&&\int _{1}^T\frac {L(iv+it_0,\chi )Q^{iv}dv}{v}\notag 
\\ &&\hspace {10mm}=\hspace {4mm}C_\chi \sqrt M\int _1^T\frac {\sqrt {|v+t_0|}}{v}
\left (\sum _{n\leq Z}\frac {\overline \chi (n)e\left (f(v)\right )}{n}\right )dv
+\mathcal O\left (\sqrt M\log Z\int _{1}^T\frac {dv}{v}\right )\notag 
\\ &&\hspace {10mm}\ll \hspace {4mm}\sqrt M\sum _{n\leq Z}\frac {1}{n}\left |
\int _{1}^TG(v)e\left (f(v)\right )dv\right |+\sqrt M\log Z\hspace {15mm}
\end {eqnarray}
subject to \eqref {goro}, where
\[ G(v)=\frac {\sqrt {|v+t_0|}}{v}.\]
We now bound the integral in \eqref {bora}.  Take $R\geq 1$.  For $v\in (R,2R)$ we 
have $|v+t_0|\ll R+|t_0|$ so from \eqref {deriv}
\begin {eqnarray}\label {1111}
v\in (R,2R)\hspace {5mm}\implies \hspace {5mm}\Bigg \{ \begin {array}{ll}G(v)&\ll \hspace {3mm}\frac {\sqrt {R+|t_0|}}{R}\hspace {5mm}\text { always}
\\ |f''(v)|&\gg \hspace {3mm}\frac {1}{R+|t_0|}\hspace {5mm}\text { if }v\not =-t_0.
\end {array}\hspace {15mm}\hspace {15mm}
\end {eqnarray}
and we now consider the various scenarios for the sizes of $R$ and $t_0$.  Suppose 
first that $R$ is large and $-t_0\in (R+1,2R-1)$.  Then the above bounds become
\begin {eqnarray*}
 v\in (R,2R)\char92 \{ t_0\} \hspace {5mm}\implies \hspace {5mm}|f''(v)|\gg \frac {1}{R}
\hspace {5mm}\text {and}\hspace {5mm}G(v)\ll \frac {1}{\sqrt R}\hspace {15mm}
\end {eqnarray*}
so that from Lemma 4.5 of \cite {titchmarsh} 
\begin {eqnarray*}
\int _R^{2R}G(v)e\left (f(v)\right )dv&=&
\left (\int _R^{-t_0-1}+\int _{-t_0-1}^{-t_0+1}+\int _{-t_0+1}^{2R}\right )
G(v)e\left (f(v)\right )dv
\\ &\ll &1
\end {eqnarray*}
having bounded the second integral crudely with \eqref {1111}.  If $-t_0\not \in (R+1,2R-1)$ then $v+t_0\not =0$ 
for $v\in (R+1,2R-1)$ so the above bounds and the same lemma imply
\begin {eqnarray*}
\int _R^{2R}G(v)e\left (f(v)\right )dv&=&\left (\int _R^{R+1}+\int _{R+1}^{2R-1}+\int _{2R-1}^{2R}\right )G(v)e\left (f(v)\right )dv
\\ &\ll &\frac {\sqrt {R+|t_0|}}{R}+\frac {\sqrt {R+|t_0|}}{R}\cdot \frac {1}{\sqrt {R+|t_0|}}+\frac {\sqrt {R+|t_0|}}{R}
\\ &\ll &1+\sqrt {|t_0|}
\end {eqnarray*}
having bounded the first and third integrals crudely with \eqref {1111}.  If $R$ is not large then 
\begin {equation}\label {yogurt3}
\int _R^{2R}G(v)e\left (f(v)\right )dv\ll 1+\sqrt {|t_0|}
\end {equation}
is clear from \eqref {1111} so we conclude that \eqref {yogurt3} holds for 
all $R\geq 1$ and subject to no constraints on $t_0$.  Consequently 
\begin {eqnarray*}
\int _1^TG(v)e\left (f(v)\right )dv\ll 1+\sqrt {|t_0|}
\end {eqnarray*}
so \eqref {bora} implies
\begin {eqnarray*}
\int _{1}^T\frac {L(iv+it_0,\chi )Q^{iv}dv}{v}&\ll &\sqrt M\Big (1+\sqrt {|t_0|}\Big )
\sum _{n\leq Z}\frac {1}{n}+\sqrt M\log Z
\\ &\ll &\sqrt M\log Z\Big (1+\sqrt {|t_0|}\Big )
\end {eqnarray*}
which proves the lemma if we set for example $Z=1+(1+T+|t_0|)^2\sqrt M$ in accordance with \eqref {goro}.
\end {proof}
Suppose $X\geq 1$ and $m\in \mathbb N$
.  For $w\in \mathbb C$ 
with $\mathfrak I\mathfrak m(w)\ll X^{\mathcal O(1)}$ 
and $\mathfrak R\mathfrak e(w)\geq 0$ and for a primitive character $\chi ^*$ modulo $m$, it is well known that
\begin {equation}\label {hh}
L(w,\chi ^*)\ll X^\epsilon \sqrt {m(1+|\mathfrak I\mathfrak m(w)|)}.
\end {equation}
\newtheorem {perron}[mvt]{Lemma}
\begin {perron}\label {perron}
For $q,d\in \mathbb N$ with $d|q$ and $s\in \mathbb C$ define
\[ U_s(q,d)=\frac {1}{q}\sum _{D|d}D^s\phi (q/D).\]
Then for any $d,M\in \mathbb N$, $t\in \mathbb R$, $x,Q\geq 1$, and any Dirichlet character $\chi $ mod $M$, we have 
for $d,M,|t|,Q\leq x^{\mathcal O(1)}$
\[ \sum _{u\leq Q}\frac {\chi \left (u/(u,d)\right )}{\left (u/(u,d)\right )^s}=\frac {U_s(dM,d)Q^{1-s}}{1-s}
+\mathcal O\left (x^\epsilon \sqrt {M(1+|t|)}\right )\]
where $s=it$, and where the main term is present if and only if $\chi $ is principal.
\\
\\ Moreover, the result remains true if $s$ is assumed to be in 
the region $\{ s\in \mathbb C|\sigma \geq 0\text { and }|s|\leq 1/2\} $.
\end {perron}
\begin {proof}
Throughout we write $s=it$ and for $w\in \mathbb C$ always $w=u+iv$, for real $u,v$.  As in the last proof we allow 
the $\ll ,\mathcal O$ symbols to contain terms up to $x^\epsilon $ (and therefore 
also $d^\epsilon ,M^\epsilon ,|t|^\epsilon ,Q^\epsilon $).
\\
\\ Let $\chi ^*$ be a primitive character of modulus $m$ say, with $m\leq x^{\mathcal O(1)}$.  Since $\chi ^*$ is principal if 
and only if $m=1$ we may define for any $X>0$
\[ R_{\chi ^*}(X)=\frac {X^{1-s}}{1-s}\left \{ \begin {array}{ll}1&\text { if }m=1
\\ 0&\text { if $\chi ^*$ is not principal}.\end {array}\right .\]
Write $A$ for the implied constant in the hypothesis and take parameters $2\leq X,T\leq x^{A+2}$ with $T$ so large that
\begin {equation}\label {oia}
T>|t|
\end {equation}
and
\begin {equation}\label {oia2}
T>X^2. 
\end {equation}
Perron's formula (Theorem 2 in Part II, Section 2 of \cite {tenen}) implies for $\kappa >1$
\begin {eqnarray}\label {reads}
\sum _{n\leq X}\frac {\chi ^*(n)}{n^{s}}&=&\frac {1}{2\pi i}\int _{\kappa \pm iT}\frac {L(w+s,\chi ^*)X^wdw}{w}
+\mathcal O\left (X^\kappa \sum _{n=1}^\infty \frac {1}{n^\kappa \big (1+T|\log (X/n)|\big )}\right )\notag 
\\ &=:&I(X,T)+\mathcal O\Big (E(X,T)\Big ).
\end {eqnarray}
If $m=1$ then $L(w+s,\chi ^*)=\zeta (w+s)$ and if $\chi ^*$ is non-principal then $L(w+s,\chi ^*)$ is holomorphic for $u>0$, so 
by the Residue Theorem and \eqref {oia}
\begin {eqnarray}\label {hv}
I(X,T)=R_{\chi ^*}(X)\hspace {3mm}-\hspace {3mm}\frac {1}{2\pi i}\left (\int _{\kappa +iT}^{iT}+\int _{\mathcal L}
+\int _{-iT}^{\kappa -iT}\right )\frac {L(w+s,\chi ^*)X^wdw}{w},\hspace {15mm}
\end {eqnarray}
where $\mathcal L$ is the vertical line from $iT$ to $-iT$ except for a 
half circle $\mathcal C$ from $\delta i$ to $-\delta i$ to the right of $0$, where $\delta =1/\log X$.  From \eqref {hh} we have
\begin {eqnarray*}
\int _{\kappa +iT}^{iT}\frac {L(w+s,\chi ^*)X^wdw}{w}&\ll &X^\kappa \int _0^\kappa \frac {|L(u+iT+it)|du}{|u+iT|}
\\ &\ll &\frac {x^\epsilon X^\kappa \sqrt {m(1+T+|t|)}}{T}
\end {eqnarray*}
and similarly for the other horizontal integral in \eqref {hv}.  For the vertical integral Lemma \ref {mvt} and \eqref {hh} imply
\begin {eqnarray*}
&&\int _{\mathcal L}\frac {L(w+s,\chi ^*)X^wdw}{w}
\\ &&\hspace {10mm}\ll \hspace {4mm}\left |\int _1^T\frac {L\left (iv+it,\chi ^*\right )X^{iv}dv}{v}\right |
+\left (\int _{\mathcal C}+\int _{\delta \leq |v|\leq 1}\right )
\frac {\left |L\left (w+s,\chi ^*\right )\right |\cdot |X^w|\cdot dw}{|w|}
\\ &&\hspace {10mm}\ll \hspace {4mm}x^\epsilon \sqrt {m(1+|t|)}+X^{\delta }\sqrt {m(1+|t|)}
\left (\int _{\mathcal C}+\int _{\delta \leq |v|\leq 1}\right )\frac {dw}{|w|}
\\ &&\hspace {10mm}\ll \hspace {4mm}x^\epsilon \sqrt {m(1+|t|)}+X^{1/\log X}\sqrt {m(1+|t|)}\cdot |\log \delta |
\\ &&\hspace {10mm}\ll \hspace {4mm}x^\epsilon \sqrt {m(1+|t|)}.
\end {eqnarray*}
Using these bounds for the integrals in \eqref {hv} and inserting the result into \eqref {reads} we get
\begin {eqnarray}\label {becomes}
\sum _{n\leq X}\frac {\chi ^*(n)}{n^{s}}&=&R_{\chi ^*}(X)+\mathcal O\left (x^\epsilon \left (\frac {X^\kappa \sqrt {m(1+T+|t|)}}{T}
+\sqrt {m(1+|t|)}+E(X,T)\right )
\right ).\hspace {15mm}
\end {eqnarray}
In general for $Z>-1$
\[ |\log (1+Z)|\geq \frac {|Z|}{1+Z}.\]
For $X/2\leq n\leq 3X/2$ we have $(n-X)/X>-1$ so that
\[ |\log (X/n)|=\left |\log \left (1+\frac {n-X}{X}\right )\right |\geq \frac {|n-X|}{n}\geq \Big \lfloor |n-X|\Big \rfloor \Big /n\]
and therefore
\begin {eqnarray}\label {tenen1}
\sum _{X/2\leq n\leq 3X/2}\frac {1}{n^\kappa |\log (X/n)|}\leq X^{1-\kappa }\left (1+2\sum _{h\leq X}\frac {1}{h}\right )
\ll X^{1-\kappa }.\hspace {15mm}
\end {eqnarray}
If $n$ is not in this range then $|\log (X/n)|\gg 1$ so from \eqref {tenen1}
\begin {eqnarray*}
X^\kappa \sum _{n=1}^\infty \frac {1}{n^\kappa \big (1+T|\log (X/n)|\big )}&\ll &
X^\kappa \left (\frac {\zeta (\kappa )}{T}+\frac {X^{1-\kappa }}{T}\right )
\\ &\ll &\frac {1}{T}\left (\frac {X^\kappa }{\kappa -1}+X\right )
\\ &\ll &\frac {X}{T}
\end {eqnarray*}
if we set $\kappa =1+1/\log X$.  Therefore $E(X,T)\ll X/T$ which we put in \eqref {becomes} to get
\begin {eqnarray*}
\sum _{n\leq X}\frac {\chi ^*(n)}{n^{s}}&=&R_{\chi ^*}(X)
+\mathcal O\left (x^\epsilon \left (\frac {X\sqrt {m(1+T+|t|)}}{T}+\sqrt {m(1+|t|)}+\frac {X}{T}\right )\right )
\\ &=&R_{\chi ^*}(X)+\mathcal O\left (x^\epsilon \sqrt {m(1+|t|)}\right )
\end {eqnarray*}
from \eqref {oia2}.  The equality obviously still valid 
if $0\leq X\leq 2$ we conclude 
that for any $0<X\leq x^{A+2}$
\begin {equation}\label {letzte}
\sum _{n\leq X}\frac {1}{n^{s}}=\frac {X^{1-s}}{1-s}
+\mathcal O\left (x^\epsilon \sqrt {(1+|t|)}\right )
\end {equation}
and
\begin {equation}\label {letzte2}
\sum _{n\leq X}\frac {\chi ^*(n)}{n^{s}}\ll x^\epsilon \sqrt {m(1+|t|)}
\end {equation}
if $\chi ^*$ is non-principal.
\\
\\ If $\chi $ is non-principal then there is an $m$ with $m|M$ and non-principal primitive character $\chi ^*$ mod $m$ for which
\[ \chi (n)=\left \{ \begin {array}{ll}\chi ^*(n)&\text { if }(n,M)=1\\ 0&\text { if not}\end {array}\right .\]
so that
\begin {eqnarray*}
\sum _{u\leq Q}\frac {\chi \left (u/(u,d)\right )}{\left (u/(u,d)\right )^s}&=&
\sum _{D|d}\sum _{u\leq Q/D\atop {(u,d/D)=1}}\frac {\chi (u)}{u^s}
\\ &=&\sum _{D|d}\sum _{u\leq Q/D\atop {(u,dM/D)=1}}\frac {\chi ^*(u)}{u^s}
\\ &=&\sum _{D|d\atop {\Delta |dM/D}}\frac {\mu (\Delta )\chi ^*(\Delta )}{\Delta ^s}\sum _{u\leq Q/D\Delta }\frac {\chi ^*(u)}{u^s}
\\ &\ll &\sum _{D,\Delta |dM}\left |\sum _{u\leq Q/D\Delta }\frac {\chi ^*(u)}{u^s}\right |
\\ &\ll &x^\epsilon \sqrt {m(1+|t|)}
\end {eqnarray*}
from \eqref {letzte2}, which proves the lemma for $\chi $ non-principal.  If $\chi $ is principal 
then we use \eqref {letzte} to deduce
\begin {eqnarray*}
\sum _{u\leq Q}\frac {\chi \left (u/(u,d)\right )}{\left (u/(u,d)\right )^s}&=&\sum _{u\leq Q\atop {\left (u/(u,d),M\right )=1}}
\frac {1}{\left (u/(u,d)\right )^s}
\\ &=&\sum _{D|d}\sum _{\Delta |dM/D}\frac {\mu (\Delta )}{\Delta ^s}\sum _{u\leq Q/D\Delta }\frac {1}{u^s}
\\ &=&\frac {Q^{1-s}}{1-s}\sum _{D|d}\frac {1}{D^{1-s}}\sum _{\Delta |dM/D}\frac {\mu (\Delta )}{\Delta }
+\mathcal O\left (x^\epsilon \sqrt {1+|t|}\right )
\\ &=&\frac {Q^{1-s}}{1-s}\sum _{D|d}\frac {1}{D^{1-s}}\cdot \frac {\phi (dM/D)}{dM/D}
+\mathcal O\left (x^\epsilon \sqrt {1+|t|}\right )
\end {eqnarray*}
which proves the lemma for $\chi $ principal.  The last claim is an easy consequence of 
partial summation and the Polya-Vinogradov Inequality.
\end {proof}

\begin {center}
\section {- \hspace {2mm} Proof of theorem}
\end {center}
Let $1\leq Q\leq x$ be given.  If $Q\leq \sqrt x$ then \eqref {vary} and the first claim of Lemma \ref {el} (ii) imply
\[ V(x,Q)\ll x^\epsilon \Big (x^{1+\theta }+x^{2\Delta }+x^{1/2+\Delta }\Big )\ll x^{1+\Delta +\epsilon }\]
which is our theorem, so we assume
\begin {equation}\label {qq}
Q>\sqrt x.
\end {equation}
Since for $(d_i^k,d_j^k)|h_i-h_j$ (as in the proof of part (iii) of Lemma \ref {et})
\[ \frac {1}{[d_1^k,...,d_r^k]}\ll \frac {1}{d_1^k\cdot \cdot \cdot d_r^k}\]
we have 
\begin {eqnarray}\label {dobby}
\sum _{d_1,...,d_r=1\atop {d_i^k-d_j^k|h_i-h_j}}^\infty \frac {|\mu (d_1)\cdot \cdot \cdot \mu (d_r)|}{[q,[d_1^k,...,d_r^k]]}
\ll \frac {1}{q}\sum _{d_1,...,d_r}\frac {(q,d_1^k\cdot \cdot \cdot d_r^k)}{d_1^k\cdot \cdot \cdot d_r^k}
\ll \frac {1}{q}\sum _{n=1}^\infty (q,n^k)n^{\epsilon -k}\ll q^{\epsilon -1}.\hspace {10mm}
\end {eqnarray}
From \eqref {vary} and \eqref {etaultra}
\begin {eqnarray}\label {radiator2}
V(x,Q)&=&\sum _{q\leq Q}\sum _{a=1}^q\sum _{n,m\leq x\atop {n,m\in \mathcal R\atop {n\equiv m\equiv a(q)}}}1
-2x\sum _{q\leq Q}\sum _{a=1}^q\eta (q,a)\sum _{n\leq x\atop {n\in \mathcal R\atop {n\equiv a(q)}}}1
+x^2\sum _{q\leq Q}\sum _{a=1}^q\eta (q,a)^2\notag 
\\ &=:&S_1(x,Q)-2xS_2(x,Q)+x^2\sum _{q\leq Q}W(q).
\end {eqnarray}
For $d_1,...,d_r\in \mathbb N$ write $d^*=[d_1,...,d_r]$.  Denote by $V$ the unique 
solution modulo $[(q,d_1),...,(q,d_r)]=(q,d^*)$ to $n\equiv -\mathbf h\left ((q,\mathbf d)\right )$.  From \eqref {etaultra} we have 
for a new paramter $X>0$
\begin {eqnarray*}
\sum _{a=1}^q\eta (q,a)\sum _{n\leq X\atop {n\in \mathcal R\atop {n\equiv a(q)}}}1=
\sum _{n\leq X\atop {n\in \mathcal R}}\eta (q,n)=\sum _{d_1,...,d_r=1\atop {d_i^k-d_j^k|h_i-h_j}}^\infty 
\frac {\mu (d_1)\cdot \cdot \cdot \mu (d_r)}{[q,{d^*}^k]}
\sum _{n\leq X\atop {n\in \mathcal R\atop {n\equiv V\left ((q,{d^*}^k)\right )}}}1.
\end {eqnarray*}
From Lemma \ref {el} (i) the inner sum here is
\[ A_{q,d^*}X+\mathcal O\left (X^{\Delta +\epsilon }\right )\]
for some $A_{q,d^*}$ and therefore from \eqref {dobby}
\begin {eqnarray}\label {ok}
\sum _{a=1}^q\eta (q,a)\sum _{n\leq X\atop {n\in \mathcal R\atop {n\equiv a(q)}}}1
&=&X\sum _{d_1,...,d_r=1\atop {d_i^k-d_j^k|h_i-h_j}}^\infty \frac {A_{q,d^*}\mu (d_1)\cdot \cdot \cdot \mu (d_r)}{[q,{d^*}^k]}
+\mathcal O\left (X^{\Delta +\epsilon }\sum _{d_1,...,d_r=1\atop {d_i^k-d_j^k|h_i-h_j}}^\infty 
\frac {|\mu (d_1)\cdot \cdot \cdot \mu (d_r)|}{[q,d^*]}\right )\notag 
\\ &=:&XB_q+\mathcal O\left (X^{\Delta +\epsilon }q^{\epsilon -1}\right ).
\end {eqnarray}
On the other hand Lemma \ref {el} (i) says
\[ \sum _{n\leq X\atop {n\in \mathcal R\atop {n\equiv a(q)}}}1=X\eta (q,a)+o(X)\]
so from \eqref {radiator2} 
\[ \sum _{a=1}^q\eta (q,a)\sum _{n\leq X\atop {n\in \mathcal R\atop {n\equiv a(q)}}}1=X\sum _{a=1}^q\eta (q,a)^2+o_q(X)
=XW(q)+o_q(X)\]
and so from \eqref {ok}
\[ XB_q+\mathcal O\left (X^{\Delta +\epsilon }q^{\epsilon -1}\right )=XW(q)+o_q(X).\]
Therefore we must have
\[ B_q=W(q)\]
and setting $X=x$ in \eqref {ok} we deduce
\begin {eqnarray*}
\sum _{a=1}^q\eta (q,a)\sum _{n\leq x\atop {n\in \mathcal R\atop {n\equiv a(q)}}}1
&=&xW(q)+\mathcal O\left (x^{\Delta +\epsilon }q^{\epsilon -1}\right )
\end {eqnarray*}
so that, from \eqref {radiator2},
\begin {eqnarray}\label {beiro}
S_2(x,Q)=\sum _{q\leq Q}\Big (xW(q)+\mathcal O\left (x^{\Delta +\epsilon }q^{\epsilon -1}\right )\Big )
=x\sum _{q\leq Q}W(q)+\mathcal O\left (x^{\Delta +\epsilon }\right ).
\end {eqnarray}
Let $\rho $ be as in Lemma \ref {bb}.  From Definition 3.3 we have $G(1)H(1,1)=1$ so 
from \eqref {exponentialsum} and Lemma \ref {el} (i)
\[ \sum _{n\leq x\atop {n\in \mathcal R}}1=\rho x+\mathcal O\left (x^{\Delta +\epsilon }\right )\]
so that 
\[ \sum _{q\leq Q}\sum _{n\leq x\atop {n\in \mathcal R}}1=\rho xQ+\mathcal O\left (x^{1+\Delta +\epsilon }\right )\]
and therefore from \eqref {radiator2}
\begin {eqnarray*}
S_1(x,Q)&=&\sum _{q\leq Q}\sum _{n,m\leq x\atop {n,m\in \mathcal R\atop {n\equiv m(q)}}}1
\\ &=&2\sum _{q\leq Q}\sum _{m<n\leq x\atop {n,m\in \mathcal R\atop {n\equiv m(q)}}}1
+\sum _{q\leq Q}\sum _{n\leq x\atop {n\in \mathcal R}}1
\\ &=:&2S_4(x,Q)+\rho xQ+\mathcal O\left (x^{1+\Delta +\epsilon }\right ).
\end {eqnarray*}
Putting this and \eqref {beiro} in \eqref {radiator2} gives
\begin {equation}\label {tan}
V(x,Q)=2S_4(x,Q)+\rho xQ-x^2\sum _{q\leq Q}W(q)+\mathcal O\left (x^{1+\Delta +\epsilon }\right )
\end {equation}
and now our task is to study $S_4(x,Q)$ using the circle method.  
\\
\\ Let $F,f$ be as in Lemma \ref {cm}.  Writing the congruence condition in $S_4(x,Q)$ out 
explicitly and using orthogonality we have
\begin {eqnarray}\label {wy}
S_4(x,Q)&=&\sum _{q\leq Q}\sum _{l\leq x/q}\sum _{n,m\leq x\atop {n,m\in \mathcal R\atop {n-m=ql}}}1
=\int _{\mathfrak U}F(-\alpha )|f(\alpha )|^2d\alpha 
\end {eqnarray}
for any unit interval $\mathfrak U$.  As 
in the comments preceeding Lemma \ref {cm}, denote by $\mathfrak F(a/q)$ the Farey arc at $a/q$ in 
the Farey dissection of order $\gamma $, where $(a,q)=1$.  Then \eqref {wy} 
and \eqref {farey999} imply
\begin {eqnarray}\label {lamp}
S_4(x,Q)&=&\sum _{q\leq 2\sqrt x}\sideset {}{'}\sum _{a=1}^q
\int _{\mathfrak F(a/q)}F(-\alpha )|f(\alpha )|^2d\alpha \notag 
\\ &&\hspace {10mm}+\hspace {4mm}\mathcal O\left (1+\sum _{2\sqrt x<q\leq \gamma }\sideset {}{'}\sum _{a=1}^q
\int _{\mathfrak F(a/q)}|F(-\alpha )|\cdot |f(\alpha )|^2d\alpha \right )\notag 
\\ &=:&M(\gamma )+\mathcal O\Big (1+E_x(\gamma )\Big ).
\end {eqnarray}
Let $\theta ,\Delta ,G,H,\rho ,I$ and $\hat J$ be as in Lemma \ref {cm} and as in that lemma write $\alpha =a/q+\beta $ 
whenever $\alpha \in \mathfrak F(a/q)$.  Suppose $2\sqrt x\leq \gamma \leq x^{3/4}$.  From part (A) of that lemma
\[ |f(\alpha )|^2=|\rho G(q)H(q,a)I(\beta )|^2+2\mathfrak R\mathfrak e\left (\overline {\rho 
G(q)H(q,a)I(\beta )}\hat J(\alpha )\right )+|\hat J(\alpha )|^2\]
so that from parts (D) and (C) we have
\begin {eqnarray}\label {pysgota}
M(\gamma )&=&|\rho |^2\sum _{q\leq 2\sqrt x}|G(q)|^2\sideset {}{'}\sum _{a=1}^q|H(q,a)|^2
\int _{\mathfrak F(a/q)}F(-\alpha )|I(\beta )|^2d\alpha  \notag 
\\ &&+\hspace {4mm}2\mathfrak R\mathfrak e\left (\overline \rho 
\sum _{q\leq 2\sqrt x}\overline {G(q)}\sideset {}{'}\sum _{a=1}^q\overline {H(q,a)}
\int _{\mathfrak F(a/q)}F(-\alpha )\overline {I(\beta )}\hat J(\alpha )d\alpha \right )\notag 
\\ &&\hspace {8mm}+\hspace {4mm}\sum _{q\leq 2\sqrt x}\sideset {}{'}\sum _{a=1}^q
\int _{\mathfrak F(a/q)}F(-\alpha )|\hat J(\alpha )|^2d\alpha \notag 
\\ &=&\rho ^2\sum _{q\leq 2\sqrt x}|G(q)|^2\sideset {}{'}\sum _{a=1}^q|H(q,a)|^2
\int _{\mathfrak F(a/q)}F(-\alpha )|I(\alpha -a/q)|^2d\alpha  
+\mathcal O\left (x^\epsilon \left (x^{1+\Delta }+\gamma x^{1/2+\theta }+\frac {x^{1+2\Delta }}{\gamma }\right )\right )\notag 
\\ &=:&M^*(\gamma )+\mathcal O\left (x^\epsilon 
\left (x^{1+\Delta }+\gamma x^{1/2+\theta }+\frac {x^{1+2\Delta }}{\gamma }\right )\right ).
\end {eqnarray}
Take a parameter $1\leq Z\leq 2\sqrt x$.  From Lemma \ref {cc}
\[ \sum _{q\leq Z}q^2|G(q)|^2\ll Z^{\theta +\epsilon }\]
so from \eqref {hbound} and Lemma \ref {cm} (F)
\begin {eqnarray}\label {sydyn}
\sum _{q\leq Z}|G(q)|^2\sideset {}{'}\sum _{a=1}^q|H(q,a)|^2
\int _{\mathfrak U(a/q)\char92 \mathfrak F(a/q)}|F(-\alpha )|\cdot |I(\alpha -a/q)|^2d\alpha 
\ll x^\epsilon \gamma ^2Z^{\theta }\hspace {10mm}
\end {eqnarray}
where $\mathfrak U(a/q)$ denotes the unit interval centered at $a/q$.  On 
the other hand Lemma \ref {cm} (B) and then orthogonality gives
\[ \sideset {}{'}\sum _{a=1}^q|H(q,a)|^2\int _{X}F(-\alpha )|I(\alpha -a/q)|^2d\alpha  
\ll x^{1+\epsilon }\int _{-1/2}^{1/2}|I(\beta )|^2d\beta \leq x^{2+\epsilon }\]
for any $X\subseteq [a/q-1/2,a/q+1/2]$, and from Lemma \ref {cc} 
\[ \sum _{q>Z}|G(q)|^2\ll Z^{\theta -2+\epsilon },\]
therefore
\begin {equation}\label {sydyn2}
\sum _{q>Z}|G(q)|^2\left |\sideset {}{'}\sum _{a=1}^q|H(q,a)|^2\int _{X}F(-\alpha )|I(\alpha -a/q)|^2d\alpha \right |
\ll x^{2+\epsilon }Z^{\theta -2}.\hspace {10mm}
\end {equation}
From \eqref {sydyn} and \eqref {sydyn2}
\begin {eqnarray}\label {gwell}
&&M^*(\gamma )=\rho ^2
\sum _{q\leq x}|G(q)|^2\sideset {}{'}\sum _{a=1}^q|H(q,a)|^2\int _{\mathfrak U(a/q)}F(-\alpha )|I(\alpha -a/q)|^2d\alpha \notag 
\\ &&\hspace {14mm}+\hspace {4mm}\mathcal O\left (\hspace {1mm}\sum _{q\leq Z}|G(q)|^2\sideset {}{'}\sum _{a=1}^q|H(q,a)|^2
\int _{\mathfrak U(a/q)\char92 \mathfrak F(a/q)}|F(-\alpha )|\cdot |I(\alpha -a/q)|^2d\alpha \right .\notag 
\\ &&\hspace {18mm}+\hspace {4mm}\sum _{q>Z}|G(q)|^2\left |\sideset {}{'}\sum _{a=1}^q|H(q,a)|^2
\int _{\mathfrak U(a/q)\char92 \mathfrak F(a/q)}F(-\alpha )|I(\alpha -a/q)|^2d\alpha \right |\notag 
\\ &&\left .\hspace {22mm}+\hspace {4mm}\sum _{q>Z}|G(q)|^2\left |\sideset {}{'}\sum _{a=1}^q|H(q,a)|^2
\int _{\mathfrak U(a/q)}F(-\alpha )|I(\alpha -a/q)|^2d\alpha \right |\hspace {1mm}\right )\notag 
\\ &&=\hspace {4mm}\rho ^2\sum _{q\leq x}|G(q)|^2\sideset {}{'}\sum _{a=1}^q|H(q,a)|^2
\int _{\mathfrak U(a/q)}F(-\alpha )|I(\alpha -a/q)|^2d\alpha 
+\mathcal O\left (x^\epsilon \left (\gamma ^2Z^{\theta }+x^2Z^{\theta -2}\right )\right )\notag 
\\ &&=\hspace {4mm}\rho ^2\sum _{q\leq x}|G(q)|^2\sideset {}{'}\sum _{a=1}^q|H(q,a)|^2
\int _{-1/2}^{1/2}F(-a/q-\beta )|I(\beta )|^2d\beta 
+\mathcal O\left (\gamma ^{2-\theta }x^{\theta +\epsilon }\right )\hspace {15mm}
\end {eqnarray}
on choosing $Z=x/\gamma $.  Since it is straightforward to establish that 
\[ \int _{-1/2}^{1/2}|I(\beta )|^2e\left (-\beta n\right )d\beta  
=
x-n+\mathcal O\left (1\right )\]
we have ($F$ is defined in Lemma \ref {cm})
\[ \int _{-1/2}^{1/2}F(-a/q-\beta )|I(\beta )|^2d\beta =\sum _{uv\leq x\atop {u\leq Q}}e(-auv/q)(x-uv)
+\mathcal O\left (x^{1+\epsilon }\right )\]
so we deduce from \eqref {gwell}, \eqref {hbound} and \eqref {gbound}
\begin {eqnarray*}
M^*(\gamma )&=&\rho ^2\sum _{q\leq x}|G(q)|^2\sideset {}{'}\sum _{a=1}^q|H(q,a)|^2
\left (\sum _{uv\leq x\atop {u\leq Q}}e(-auv/q)\Big (x-uv\Big )+\mathcal O(x^{1+\epsilon })\right )
+\mathcal O\left (\gamma ^{2-\theta }x^{\theta +\epsilon }\right )
\\ &=&\rho ^2\sum _{q\leq x}|G(q)|^2\sum _{uv\leq x\atop {u\leq Q}}\Big (x-uv\Big )\Phi _q(-uv)
+\mathcal O\left (x^\epsilon \Big (x+\gamma ^{2-\theta }x^{\theta }\Big )\right ),
\end {eqnarray*}
where $\Phi _q(n)$ is as in Lemma \ref {bb}, so from \eqref {lamp} and \eqref {pysgota} 
\begin {eqnarray}\label {na}
S_4(x,Q)&=&\rho ^2\sum _{q\leq x}|G(q)|^2\sum _{uv\leq x\atop {u\leq Q}}\Big (x-uv\Big )\Phi _q(-uv)\notag 
\\ &&\hspace {10mm}+\hspace {4mm}\mathcal O\left (x^\epsilon \left (x^{1+\Delta }+\gamma x^{1/2+\theta }
+\frac {x^{1+2\Delta }}{\gamma }+E_x(\gamma )+\gamma ^{2-\theta }x^{\theta }\right )\right ).\hspace {15mm}
\end {eqnarray}
Recall that $\theta =1/k$ and $\Delta =2/(k+1)$.  If $k>2$ we 
set $\gamma =2\sqrt x$ so that $E_x(\gamma )=0$ and $\gamma \geq x^\Delta $ to deduce that
\begin {eqnarray}\label {dymani}
S_4(x,Q)&=&\rho ^2\sum _{q\leq x}|G(q)|^2\sum _{uv\leq x\atop {u\leq Q}}\Big (x-uv\Big )\Phi _q(-uv)
+\mathcal O\left (x^\epsilon \Big (x^{1+\Delta }+x^{1+\theta }\Big )\right )\notag 
\\ &=:&\rho ^2\mathcal J(x,Q)+\mathcal O\left (x^{1+\Delta +\epsilon }\right ).
\end {eqnarray}
If $k=2$ we set $\gamma =x^{2/3}$ and deduce from Lemma \ref {cm} (E) that the error term in \eqref {na} is up to 
an $x^\epsilon $ bound
\[ \ll x^{5/3}+\gamma x+\frac {x^{7/3}}{\gamma }+x\gamma +\gamma ^{3/2}x^{1/2}\ll x^{5/3}=x^{1+\Delta }\]
to conclude that \eqref {dymani} holds for all $k\geq 2$.  
\\
\\ This finishes our circle method work and it 
remains to evaluate $\mathcal J(x,Q)$.  We use the periodicity of $\Phi _q(n)$ modulo $q$ and apply Perron's formula to 
evaluate precisely the remaining quantity.
\\
\\ We 
make the convention that whenever we have the letter $\mathcal D$ appearing in a context involving natural numbers $q,a$ we 
mean $\mathcal D=(q,a)$.  For any $u\in \mathbb N$ we then write $u'=u/(u,\mathcal D)$.  Sorting 
the $uv$ according to the residue $a$ mod $q$ we have 
\begin {eqnarray}\label {uno}
\mathcal J(x,Q)&=&\sum _{q\leq x}|G(q)|^2\sum _{a=1}^q\Phi _q(-a)\sum _{uv\leq x\atop {u\leq Q\atop {uv\equiv a(q)}}}
\Big (x-uv\Big )\notag 
\\ &=&\sum _{q\leq x}|G(q)|^2\sum _{a=1}^{q}\Phi _q(-a)
\sum _{uv\leq x\atop {u\leq Q\atop {\mathcal D|uv\atop {uv/\mathcal D\equiv a'(q')}}}}\Big (x-uv\Big )\notag 
\\ &=:&\sum _{q\leq x}|G(q)|^2\sum _{a=1}^{q}\Phi _q(-a)\mathcal U(q,a).
\end {eqnarray}
For $n\in \mathbb N$ with $(n,q')=1$ denote by $\overline n$ the inverse of $n$ modulo $q'$.  We have
\begin {eqnarray}\label {sebon}
\mathcal U(q,a)&=&\sum _{u\leq Q}
\sum _{v\leq x/u\atop {\mathcal D/(\mathcal D,u)|v\atop {uv/\mathcal D\equiv a'(q')}}}\Big (x-uv\Big )\notag 
\\ &=&\sum _{u\leq Q\atop {(u',q')=1}}\mathcal Du'\sum _{v\leq x/\mathcal Du'\atop {v\equiv \overline {u'}a'(q')}}\Big (x/\mathcal Du'-v\Big )\notag 
\\ &=:&\sum _{u\leq Q\atop {(u',q')=1}}\mathcal Du'\mathcal V_{q,a}(u).
\end {eqnarray}
Through the orthogonality of Dirichlet characters and a Perron 
formula (taking $w=1$ in (11) of Section 2, Part II in \cite {tenen}, page 134, the 
relevant quantities being defined at the start of that section) we have 
\begin {eqnarray}\label {sebon2}
\mathcal V_{q,a}(u)&=&\frac {1}{\phi (q')}
\sum _{\chi }\overline \chi (\overline {u'}a')\sum _{v\leq x/\mathcal Du'}\chi (v)\Big (x/\mathcal Du'-v\Big )\notag 
\\ &=&\frac {1}{\phi (q')}\sum _{\chi }\overline \chi (\overline {u'}a')
\int _{2-i\infty }^{2+i\infty }\frac {L(s,\chi )}{s(s+1)}\left (\frac {x}{\mathcal Du'}\right )^{s+1}ds;
\end {eqnarray}
here and in what follows the sum $\Sigma _\chi $ runs over the Dirichlet characters modulo $q'$ and for $s\in \mathbb C$ we always 
write $s=\sigma +it$ for real numbers $\sigma ,t$.  Denote 
by $\mathcal L$ the contour from $-ix^6$ to $ix^6$ which is a vertical line except for a small 
detour to the right of $0$.  Define
\[ \delta (\chi )=\left \{ \begin {array}{ll}
1&\text { if }\chi =\chi _0\\ 0&\text { if }\chi \not =\chi _0.\end {array}\right .\]
The 
parts of the above integral with $|t|\geq x^6$ contribute to the integral
\[ \ll x^3\int _{|t|\geq x^6}\frac {|L(2+it,\chi )|dt}{t^2}\ll x^3\int _{|t|\geq x^6}\frac {dt}{t^2}\ll \frac {1}{x^3}\]
so pulling the remaining part of the integral to the left, and so picking up a simple pole at $s=1$ if $\chi =\chi _0$, we see that the 
integral in \eqref {sebon2} is
\begin {eqnarray*}
&&\int _{\mathcal L}\frac {L(s,\chi )}{s(s+1)}\left (\frac {x}{\mathcal Du'}\right )^{s+1}ds
+\frac {Res_{s=1}L(s,\chi _0)}{2}\left (\frac {x}{\mathcal Du'}\right )^2\delta (\chi )+\mathcal O\left (\frac {1}{x^3}\right )
\\ &&\hspace {10mm}=\hspace {4mm}\int _{\mathcal L}\frac {L(s,\chi )}{s(s+1)}\left (\frac {x}{\mathcal Du'}\right )^{s+1}ds
+\frac {\phi (q')}{2q'}\left (\frac {x}{\mathcal Du'}\right )^2\delta (\chi )+\mathcal O\left (\frac {1}{x^3}\right )
\end {eqnarray*}
so that for $(u',q')=1$
\begin {eqnarray*}
\mathcal V_{q,a}(u)
&=&\frac {1}{\phi (q')}\sum _{\chi }\chi (u')\overline \chi (a')
\int _{\mathcal L}\frac {L(s,\chi )}{s(s+1)}\left (\frac {x}{\mathcal Du'}\right )^{s+1}ds
\\ &&\hspace {10mm}+\hspace {4mm}\frac {1}{2q'}\left (\frac {x}{\mathcal Du'}\right )^2+\mathcal O\left (\frac {1}{x^3\phi (q')}\sum _{\chi }1\right )
\end {eqnarray*}
and so from \eqref {sebon} for $q\leq x$
\begin {eqnarray}\label {deffy}
\mathcal U(q,a)&=&\frac {1}{\phi (q')}\int _{\mathcal L}\frac {x^{s+1}}{s(s+1)\mathcal D^s}\left (
\sum _{\chi }\overline \chi (a')L(s,\chi )
\sum _{u\leq Q}\frac {\chi (u')}{{u'}^s}\right )ds\notag 
\\ &&\hspace {10mm}+\hspace {4mm}\frac {x^2}{2q'}\sum _{u\leq Q\atop {(u',q')=1}}\frac {1}{\mathcal Du'}
+\mathcal O\left (\frac {1}{x^3}\sum _{u\leq Q}\mathcal Du'\right )\notag 
\\ &=:&\mathcal I(q,a)
+\frac {x^2}{2q}\sum _{u\leq Q\atop {(u/(u,\mathcal D),q/\mathcal D)=1}}\frac {(u,\mathcal D)}{u}+\mathcal O\left (1\right ).
\end {eqnarray}
We have
\begin {equation}\label {lori}
L(s,\chi _0)=\zeta (s)\prod _{p|q'}(1-p^{-s})=:\zeta (s)\omega _s(q');
\end {equation}
define for $d|q$
\begin {equation}\label {hawsach}
\theta _s(q,d)=\frac {\omega _s(q/d)U_s(q,d)}{d^s\phi (q/d)}
\end {equation}
where $U_s(q,d)$ is as in Lemma \ref {perron}.  For $q',|t|\leq x$ we have the standard estimate 
\[ L(s,\chi )\ll x^\epsilon \sqrt {q'(1+|t|)},\hspace {10mm}0\leq \sigma \leq 1,\]
so that with Lemma \ref {perron} we see that the term in the brackets in $\mathcal I(q,a)$ is for $q\leq x$  
\begin {eqnarray}\label {wanhwn}
&=&\frac {\overline \chi _0(a')L(s,\chi _0)U_s(\mathcal Dq',\mathcal D)Q^{1-s}}{1-s}
+\mathcal O\left (x^\epsilon \sqrt {q'(1+|t|)}\sum _{\chi }|\overline \chi (a')L(s,\chi )|\right )\notag 
\\ &=&\frac {Q^2\zeta (s)\mathcal D^s\phi (q')\theta _s(q,\mathcal D)}{(1-s)Q^{s+1}}+\mathcal O\Big (x^\epsilon q'(1+|t|)\phi (q')\Big )
\end {eqnarray}
and so 
\begin {eqnarray*}
\mathcal I(q,a)&=&Q^2\int _{\mathcal L}\frac {\zeta (s)\theta _s(q,\mathcal D)}{s(s+1)(1-s)}\left (\frac {x}{Q}\right )^{s+1}ds
+\mathcal O\left (x^\epsilon q'\int _{\mathcal L}\frac {\big |(1+|t|)x^{s+1}\big |ds}{|s(s+1)\mathcal D^s|}\right )
\\ &=&-Q^2\int _{\mathcal L}\frac {\zeta (s)\theta _s(q,\mathcal D)}{(s-1)s(s+1)}
\left (\frac {x}{Q}\right )^{s+1}ds+\mathcal O\left (x^{1+\epsilon }q\right )
\end {eqnarray*}
so from Lemma \ref {bb} (iii) we have for $q\leq x$
\begin {eqnarray*}
\sum _{a=1}^{q}\Phi _q(-a)\mathcal I(q,a)
&=&-Q^2\int _{\mathcal L}\frac {\zeta (s)(x/Q)^{s+1}}{(s-1)s(s+1)}
\left (\sum _{a=1}^q\theta _s(q,\mathcal D)\Phi _q(-a)\right )ds
+\mathcal O\left (x^{1+\epsilon }q^2\right )
\end {eqnarray*}
and therefore from Lemma \ref {cc} (and since \eqref {gbound} says $|G(q)|\ll $
\begin {eqnarray}\label {aii}
&&\sum _{q\leq x}|G(q)|^2\sum _{a=1}^{q}\Phi _q(-a)\mathcal I(q,a)\notag 
\\ &&\hspace {20mm}=\hspace {4mm}-Q^2\int _{\mathcal L}\frac {\zeta (s)(x/Q)^{s+1}}{(s-1)s(s+1)}
\left (\sum _{q\leq x}|G(q)|^2\sum _{a=1}^q\theta _s(q,\mathcal D)\Phi _q(-a)\right )
ds\notag 
\\ &&\hspace {36mm}+\hspace {4mm}\mathcal O\left (x^{1+\epsilon }\sum _{q\leq x}q^2|G(q)|^2\right )\notag 
\\ &&\hspace {20mm}=:\hspace {4mm}-Q^2\int _{\mathcal L}\frac {\zeta (s)(x/Q)^{s+1}}{(s-1)s(s+1)}
\left (\sum _{q\leq x}|G(q)|^2\Delta _s(q)\right )ds
+\mathcal O\left (x^{1+\theta +\epsilon }\right ).\hspace {15mm}
\end {eqnarray}
For $\sigma \geq 0$ and $d|q$ it is clear that $U_s(q,d)\ll |d^s|q^\epsilon $ and 
so from \eqref {hawsach} that $\theta _s(q,d)\ll q^\epsilon $, so from Lemma \ref {bb} (iii) we have
\[ \Delta _s(q)\ll q^{1+\epsilon }\]
and therefore from Lemma \ref {cc}
\[ \sum _{q>x}|G(q)|^2\Delta _s(q)\ll x^{\theta -1+\epsilon }\]
so we can add in these terms to \eqref {aii} at the cost of an error of size
\[ \ll x^{\theta -1+\epsilon }Q^2\int _{\mathcal L}\left |\frac {\zeta (s)(x/Q)^{s+1}}{(s-1)s(s+1)}\right |ds
\ll x^{\theta +\epsilon }Q\int _{\pm \infty }\frac {|t|^{1/2}dt}{1+|t|^3}\ll x^{1+\theta +\epsilon }\]
to get
\begin {eqnarray}\label {VVV}
\sum _{q\leq x}|G(q)|^2\sum _{a=1}^{q}\Phi _q(-a)
\mathcal I(q,a)&=&-Q^2\int _{\mathcal L}\frac {\zeta (s)(x/Q)^{s+1}}{(s-1)s(s+1)}
\left (\sum _{q=1}^\infty |G(q)|^2\Delta _s(q)\right )ds
+\mathcal O\left (x^{1+\theta +\epsilon }\right )\notag 
\\ &=:&-Q^2\int _{\mathcal L}\frac {\zeta (s)(x/Q)^{s+1}\mathcal G(s)ds}{(s-1)s(s+1)}
+\mathcal O\left (x^{1+\theta +\epsilon }\right )\notag 
\\ &=:&-Q^2\mathcal V(x/Q)+\mathcal O\left (x^{1+\theta +\epsilon }\right ),
\end {eqnarray}
where $\mathcal G(s)$ converges absolutely (at least) for $\sigma \geq 0$.  The last 
equality with \eqref {uno} and \eqref {deffy} implies
\begin {eqnarray}\label {WWW}
\mathcal J(x,Q)&=&-Q^2\mathcal V(x/Q)
+\frac {x^2}{2}\sum _{q\leq x}\frac {|G(q)|^2}{q}\sum _{a=1}^q\Phi _q(-a)\sum _{u\leq Q\atop {(u/(u,\mathcal D),q/\mathcal D)=1}}\frac {(u,\mathcal D)}{u}\notag 
\\ &&\hspace {10mm}+\hspace {4mm}\mathcal O\left (x^{1+\theta +\epsilon }
+x^{1+\theta +\epsilon }\sum _{q\leq x}|G(q)|^2\sum _{a=1}^q|\Phi _q(-a)|\right )\notag 
\\ &=:&-Q^2\mathcal V(x/Q)+\frac {x^2\mathcal W(Q)}{2}+\mathcal O\left (x^{1+\theta +\epsilon }\right )
\end {eqnarray}
with Lemma \ref {bb} (iii) and Lemma \ref {cc}.  Our application of Perron's formula is complete and now 
our task is now to evaluate $\mathcal G(s)$.  The main point is we can write down 
an analytic continuation for this thanks to the explicit expressions for the Gauss sum in Section 3.  Recall 
the assumption of our theorem: we 
always have $R_p<p^k$, where $R_p$ is the number of distinct residue classes represented by the $h_1,...,h_r$.  
\\
\\ Recall from Lemma \ref {perron} that
\begin {equation}\label {UUU}
U_s(q,d)=\frac {1}{q}\sum _{D|d}D^s\phi (q/D)=:\frac {F^*_s(q,d)}{q}.
\end {equation}
It is straightforward to establish that $U_s(q,d)$ 
satisfies 
$U_s(qq',dd')=U_s(q,d)U_s(q',d')$ for $(q,q')=1$ and $\mathbf d|\mathbf q$ so from \eqref {hawsach} the same must be true 
of $\theta _s\left (q,d\right )$, 
and therefore Lemma \ref {bb} (i) says 
that $\Delta _s(q)$ is multiplicative (this is defined in \eqref {aii}), so from Definition \ref {qm} we have for $\sigma \geq 0$ 
\begin {eqnarray}\label {888}
\mathcal G(s)&=&\prod _p\left (\sum _{t\geq 0}|G(q)|^2\Delta _s(q)\right )
=\prod _p\left (1+\frac {1}{p^{2k}(1-R_p/p^k)^2}\sum _{1\not =q|p^k}\Delta _s(q)\right ).\hspace {15mm}
\end {eqnarray}
Define $\rho $ as in Lemma \ref {bb} (iv), namely
\[ \rho =\prod _p\left (1-\frac {R_p}{p^k}\right )
;\]
then from \eqref {888} for $\sigma \geq 0$
\begin {eqnarray}\label {dss}
\rho ^2\mathcal G(s)&=&\prod _p\left (\left (1-\frac {R_p}{p^k}\right )^2
+\frac {1}{p^{2k}}\sum _{1\not =q|p^k}\Delta _s(q)\right ).\hspace {10mm}
\end {eqnarray}
From \eqref {hawsach} and \eqref {UUU} we have
\[  \theta _s\left (q,(q,a)\right )=\frac {\omega _s(q/\mathcal D)F^*_s(q,\mathcal D)}{q\mathcal D^s\phi (q/\mathcal D)}\]
so from \eqref {aii}
\begin {eqnarray}\label {3ogloch}
\Delta _s(q)
&=&\frac {1}{q}\sum _{d|q}\frac {\omega _s(q/d)F^*_s(q,d)}{d^s\phi (q/d)}\sideset {}{'}\sum _{a=1}^{q/d}\Phi _q(-ad).
\end {eqnarray}
Define $H(q,a)$ and $\Phi (q)$ as in Lemma \ref {bb}, and take 
a prime $p$.  Denote the different residues represented by $h_1,...,h_r$ modulo $p^k$ by $H_1,...,H_{R_p}$.  For $q|p^k$ 
\[ H(q,a)=\sum _{n=1}^{R_p}e\left (\frac {aH_n}{q}\right )\]
so that
\[ \Phi (q)=\sum _{n,n'=1}^{R_p}c_q(H_n-H_{n'})\]
so from Lemma \ref {bb} (ii) 
\[ \sideset {}{'}\sum _{a=1}^{q/d}\Phi _q(-ad)=\mu (q/d)\sum _{n,n'=1}^{R_p}c_q(H_n-H_{n'})\]
and therefore from \eqref {3ogloch} 
\begin {eqnarray}\label {batri}
\Delta _s(q)&=&\frac {1}{q}\left (\sum _{d|q}\frac {\mu (q/d)\omega _s(q/d)F^*(q,d)}{d^s\phi (q/d)}\right )
\left (\sum _{n,n'=1}^{R_p}c_q(H_n-H_{n'})\right )\notag 
\\ &=:&\frac {P_s(q)}{q}\sum _{n,n'=1}^{R_p}c_q(H_n-H_{n'}).
\end {eqnarray}
Simple calculations show
\[ F_s^*(q,q)=q^s+\frac {\phi (q)\left (q^{s-1}-1\right )}{p^{s-1}-1}\]
and
\[ F_s^*(q,q/p)=F_s^*(q,q)-q^s\]
so that
\begin {eqnarray*}
P_s(q)&=&\frac {F^*_s(q,q)}{q^s}-\frac {(1-p^{-s})F^*_s(q,q/p)}{(q/p)^s\phi (p)}
\\ &=&\left (q^s+\frac {\phi (q)\left (q^{s-1}-1\right )}{p^{s-1}-1}\right )\left (\frac {1}{q^s}
-\frac {1-p^{-s}}{(q/p)^s\phi (p)}\right )+\frac {q^s(1-p^{-s})}{(q/p)^s\phi (p)}
\\ &=&q^{1-s}
\end {eqnarray*}
so from \eqref {batri} 
\[ \Delta _s(q)=q^{-s}\sum _{n,n'=1}^{R_p}c_q(H_n-H_{n'})\]
for $q|p^k$.  
The sum here is
\begin {eqnarray}\label {franzi}
\sum _{n,n'=1\atop {n\not =n'}}^{R_p}c_q(H_n-H_{n'})
+R_p\phi (q)=:\Delta ^*(q)+R_p\phi (q)
\end {eqnarray}
so 
\[ \Delta _s(q)=q^{-s}\Big (\Delta ^*(q)+R_p\phi (q)\Big )\]
and therefore, writing $X=p$ and $Y=p^{-s}$,
\begin {eqnarray}\label {franzi2}
\sum _{1\not =q|p^k}\Delta _s(q)&=&\sum _{1\not =q|p^k}q^{-s}\Delta ^*(q)
+R_p\left (\sum _{1\not =q|p^{k-1}}q^{-s}\phi (q)+p^{-sk}\phi (p^k)\right )\notag 
\\ &=&\sum _{1\not =q|p^k}q^{-s}\Delta ^*(q)+
R_p\left (\frac {(1-1/p)(p^{k(1-s)}-p^{1-s})}{p^{1-s}-1}-p^{k(1-s)-1}\right )+R_pp^{k(1-s)}\notag 
\\ &=&\sum _{1\leq t\leq k}Y^t\Delta ^*(X^t)+R_p\left (\frac {(1-1/X)((XY)^k-XY)}{XY-1}-\frac {(XY)^k}{X}\right )+R_p(XY)^k\notag 
\\ &=:&P_1(X,Y)+R_p(XY)^k
\end {eqnarray}
so from \eqref {dss} for $\sigma \geq 0$
\begin {eqnarray}\label {dsss}
\rho ^2\mathcal G(s)&=&\prod _p\left (1-\frac {2R_p}{X^k}+\frac {R_p^2}{X^{2k}}
+\frac {P_1(X,Y)}{X^{2k}}+R_p\left (\frac {Y}{X}\right )^k\right )\notag 
\\ &=&\prod _p\left (\left (1+(Y/X)^k\right )^{R_p}+1+R_p\left (\frac {Y}{X}\right )^k-\left (1+(Y/X)^k\right )^{R_p}
-\frac {2R_p}{X^k}+\frac {R_p^2}{X^{2k}}+\frac {P_1(X,Y)}{X^{2k}}\right )\notag 
\\ &=:&\prod _p\left (\left (1+(Y/X)^k\right )^{R_p}+P_2(X,Y)+\frac {P_1(X,Y)}{X^{2k}}\right );
\end {eqnarray} 
in particular
\begin {equation}\label {ana}
\rho ^2\mathcal G(0)=\prod _p\left (1-\frac {2R_p}{X^k}+\frac {R_p^2}{X^{2k}}
+\frac {P_1(X,0)}{X^{2k}}+\frac {R_p}{X^k}\right ).\hspace {10mm}
\end {equation}
For $q$ a power of $p$ and $n\in \mathbb N$ with $q\not |\hspace {1mm}n$ we 
have $c_q(n)\ll q/p$, so for $q|p^k$ all the summands in $\Delta ^*(q)$ are $\ll p^{k-1}$.  Therefore 
\begin {equation}\label {itcrowd}
\sum _{1\leq t\leq k}Y^t\Delta ^*(X^t)\ll R_p^2X^{k-1}\sum _{1\leq t\leq k}|Y|^t\ll X^{k-1}\Big (|Y|+|Y|^{k}\Big ).
\end {equation}
Fix $0<\delta <1/2k$.  If $\sigma \in [-1,1/2]$ 
then $|Y|/X\leq 1$.  
Moreover $5X|Y|/7\geq 5\sqrt 2/7>1$ so with \eqref {itcrowd}
\begin {eqnarray}\label {aufjedenFall}
P_1(X,Y)&\ll &X^{k-1}\Big (|Y|+|Y|^{k}\Big )+\underbrace {\frac {X^{k}|Y|^{k}+X|Y|}{X|Y|}}_{\ll X^{k-1}(|Y|+|Y|^k)+1}
+\frac {X^k|Y|^k}{X}\notag 
\\ &\ll &X^k+X^{2k-1}\left (\frac {|Y|}{X}\right )^k\notag 
\\ &\ll &\frac {X^{2k}}{X^{1+\delta k}},\hspace {10mm}\text {for }\sigma \geq -1+\delta .
\end {eqnarray}
From the Binomial Theorem (and even if $R_p=1$)
\begin {eqnarray}\label {hwnhefyd}
\left (1+(Y/X)^k\right )^{R_p}&=&1+R_p\left (\frac {Y}{X}\right )^k
+\mathcal O\left (\left (\frac {Y}{X}\right )^{2k}+\cdot \cdot \cdot +\left (\frac {Y}{X}\right )^{R_pk}\right )\notag 
\\ &=&1+R_p\left (\frac {Y}{X}\right )^k+\mathcal O\left (\frac {1}{p^{1+2k\delta }}\right ),
\hspace {7.5mm}\text {for }\sigma \geq -1+\frac {1}{2k}+\delta ,
\end {eqnarray}
so 
\[ P_2(X,Y)\ll \frac {1}{p^{1+2k\delta }},
\hspace {10mm}\text {for }\sigma \geq -1+\frac {1}{2k}+\delta ,\]
so that with \eqref {aufjedenFall} we have
\begin {equation}\label {r2}
P_2(X,Y)+\frac {P_1(X,Y)}{X^{2k}}\ll \frac {1}{p^{1+\delta }},\hspace {10mm}\text {for }\sigma \geq -1+\frac {1}{2k}+\delta .
\end {equation}
For $\sigma \geq -1+
\delta $ we have $|(Y/X)^k|\leq 1/p^{
\delta k}$ so \eqref {hwnhefyd} says $\left (1+(Y/X)^k\right )^{R_p}\gg 1$ 
so we deduce from \eqref {dsss} and \eqref {r2} that for $\sigma \geq 0$
\begin {eqnarray*}
\rho ^2\mathcal G(s)&=&\prod _p\Big (1+(Y/X)^k\Big )^{R_p}
\prod _p\left (1+\frac {P_2(X,Y)+P_1(X,Y)/X^{2k}}{(1+(Y/X)^k)^{R_p}}\right )
\\ &=&\prod _p\Big (1+(Y/X)^k\Big )^{r}\prod _{p^k\leq h_r}\frac {\left (1+(Y/X)^k\right )^{R_p}}{\left (1+(Y/X)^k\right )^{r}}
\prod _p\left (1+\frac {P_2(X,Y)+P_1(X,Y)/X^{2k}}{(1+(Y/X)^k)^{R_p}}\right )
\\ &=:&\frac {\zeta (sk+k)^r\mathcal F(s)}{\zeta (2sk+2k)^r}
\end {eqnarray*}
where the third product in the second line is absolutely convergent and uniformly bounded for $\sigma \geq -1+1/2k+\delta $, so 
that $\mathcal F(s)$ is holomorphic and uniformly bounded for $\sigma \geq -1+1/2k+\delta $.  
\\
\\ We now have our analytic extension for $\mathcal G(s)$ in place.  With this extension and in 
view of the clear fact 
\[ \frac {\zeta (s)\zeta (sk+k)^r}{\zeta (2sk+2k)^r}\ll |t|^{3/2+\epsilon }\log ^r|t|,
\hspace {5mm}\text { for }\sigma \geq -1+1/k\text { and }|t|\geq 1,\]
we see that the integral in $V(y)$ (see \eqref {VVV}) is 
certainly absolutely convergent for $\sigma \geq -1+1/k+\delta $ and so we may pull it to the left, picking up a simple 
pole at $s=0$, to deduce
\begin {eqnarray*}
\rho ^2\mathcal V(y)&=&-\rho ^2\zeta (0)\mathcal G(0)y+\int _{(-1+1/k+\delta )}
\frac {\zeta (s)\zeta (sk+k)^r\mathcal F(s)y^{s+1}ds}{(s-1)s(s+1)\zeta (2sk+2k)^r}
\\ &=&\frac {\rho ^2\mathcal G(0)y}{2}+\frac {1}{k}\int _{(1+k\delta )}\frac {\zeta (-1+s/k)\zeta (s)^r\mathcal F(-1+s/k)y^{s/k}ds}
{(-2+s/k)(-1+s/k)(s/k)\zeta (2s)^r}
\\ &=:&\frac {\rho ^2\mathcal G(0)y}{2}+\int _{(1+k\delta )}f(s)y^{s/k}ds.
\end {eqnarray*}
Since for $1/2\leq \sigma \leq 3/2$ and $|t|\geq 1$ we have the standard bounds $\zeta (-1+s/k)\ll |t|^{3/2-\sigma /k+\epsilon }$ 
and $\zeta (2s)\gg 1/(\log |t|)^{7}$ we see that
\[ f(s)\ll \frac {(1+|t|^{3/2-\sigma /k+\epsilon })|\zeta (s)|^r}{(1+|t|)^3}\ll \frac {|\zeta (s)|^r}{(1+|t|)^{3/2}},
\hspace {10mm}\text {for }\hspace {5mm}\frac {1}{2}\leq \sigma \leq \frac {3}{2}
.\]
By the definition of $\mathfrak c$ the integral above therefore converges absolutely for $\sigma \geq \mathfrak c$, and 
we may move the line of integration to $\sigma =\mathfrak c$, picking 
up a pole at $s=1$, to deduce 
\begin {eqnarray}\label {suarez}
\rho ^2\mathcal V(y)&=&\frac {\rho ^2\mathcal G(0)y}{2}+Res_{s=1}\left (f(s)y^{s/k}\right )+\int _{(\mathfrak c)}f(s)y^{s/k}ds\notag 
\\ &=&\frac {\rho ^2\mathcal G(0)y}{2}+Res_{s=1}\left (f(s)y^{s/k}\right )+\mathcal O_\mathfrak c\left (y^{\mathfrak c/k}\right ).
\end {eqnarray}
Since the pole of $f$ is of order $r$ a standard formula from complex analysis tells us that
\begin {eqnarray*}
Res_{s=1}\left (f(s)y^{s/k}\right )&=&\frac {1}{(r-1)!}\sum _{i+j=r-1}\left (\left (\frac {d}{ds}\right )^{(j)}
(y^{s/k})\Big |_{s=1}\right )\Bigg (Res _{s=1}\Big ((s-1)^{r-i-1}f(s)\Big )\Bigg )
\\ &=&\frac {y^{1/k}}{(r-1)!}\sum _{i+j=r-1}\left (\frac {\log y}{k}\right )^j\Bigg (Res _{s=1}\Big ((s-1)^{r-i-1}f(s)\Big )\Bigg )
\\ &=&-\frac {y^{1/k}P(\log y)}{2},
\end {eqnarray*}
where $P=P_r$ is a polynomial of degree at most $r-1$, so from \eqref {suarez}
\begin {equation}\label {muse}
-2\rho ^2\mathcal V(y)=-\rho ^2\mathcal G(0)y+y^{1/k}P(\log y)+\mathcal O\left (y^{\mathfrak c/k}\right ).\hspace {10mm}
\end {equation}
It is straightforward to establish that for $q\not |\hspace {1mm}n$
\[ \sum _{d|q}c_d(n)=0\]
so from \eqref {franzi} 
\[ \sum _{1\leq t\leq k}\Delta ^*(p^t)=\sum _{n,n'=1\atop {n\not =n'}}^{R_p}
\sum _{d|p^k\atop {d\not =1}}c_d(H_n-H_{n'})=-R_p(R_p-1)\]
so from \eqref {franzi2}
\[ P_1(X,0)=\sum _{1\leq t\leq k}\Delta ^*(p^t)+R_p\left (\frac {(1-1/X)(X^k-X)}{X-1}-X^{k-1}\right )=-R_p^2\]
and so from \eqref {ana}
\begin {equation}\label {llianbwr}
\rho ^2\mathcal G(0)=\prod _p\left (1-\frac {R_p}{X^k}\right )=\rho .
\end {equation}
From \eqref {WWW}
\begin {eqnarray}\label {W}
\mathcal W(Q)&=&\sum _{u\leq Q}\frac {1}{u}\sum _{q\leq x}\frac {|G(q)|^2}{q}
\sum _{a=1\atop {(u/(u,\mathcal D),q/\mathcal D)=1}}^{q}(u,\mathcal D)\Phi _q(-a)\notag 
\\ &=&\sum _{u\leq Q}\frac {1}{u}\sum _{q\leq x}\frac {|G(q)|^2(u,q)}{q}
\sum _{a=1\atop {(u,q)=(u,\mathcal D)}}^{q}\Phi _q(-a)\notag 
\\ &=:&\sum _{u\leq Q}\frac {1}{u}\sum _{q\leq x}\frac {|G(q)|^2f_u(q)}{q}.
\end {eqnarray}
From Lemma \ref {bb} (ii) we see that $f_u(q)$ is multiplicative and for prime powers $q$ Lemma \ref {bb} (ii) implies
\begin {eqnarray*}
f_u(q)&=&\sum _{d|q\atop {(u,q)=(u,d)}}\sideset {}{'}\sum _{a=1}^{q/d}\Phi _q(-ad)
\\ &=&\Phi (q)\sum _{d|q\atop {(u,q)=(u,d)}}\mu (q/d)
\\ &=:&(u,q)\Phi (q)F_u(q)
\end {eqnarray*}
so that for general $q$
\begin {eqnarray*}
f_u(q)&=&(u,q)\Phi (q)\left (\prod _{p^\beta ||q}F_u(p^\beta )\right ).
\end {eqnarray*}
If $p^\beta |u$ then the summation condition in $F_u(p^\beta )$ is impossible unless $d=q$ and 
so $F_u(p^\beta )=1$.  If $p^\beta \not |\hspace {1mm}u$ then the condition holds for $d=q$ and $d=q/p$ so $F_u(p^\beta )=0$.  Therefore
\[ f_u(q)=\left \{ \begin {array}{ll}q\Phi (q)&\text { if }q|u
\\ 0&\text { if not}\end {array}\right .\]
and so from \eqref {W}
\begin {eqnarray}\label {primes}
\mathcal W(Q)=\sum _{u\leq Q}\frac {1}{u}\sum _{q|u}\Phi (q)|G(q)|^2=\frac {1}{\rho ^2}\sum _{u\leq Q}W(u)\hspace {10mm}
\end {eqnarray}
from Lemma \ref {bb} (iv) and \eqref {radiator2}.  Our theorem \ref {walcott} now 
follows from \eqref {tan}, \eqref {dymani}, \eqref {WWW}, \eqref {muse}, \eqref {llianbwr} 
and \eqref {primes}.
\\
\\ 

$\hspace {1mm}$
\\
\\
\\
\\ \emph {e-mail address} - tomos.parry1729@hotmail.co.uk

\end {document}